\documentclass[a4paper, 11pt]{article} 
\usepackage{amsfonts,amsmath,amssymb,amscd}
\usepackage{latexsym}
\usepackage{graphicx}
\usepackage{color} 

\newcommand{\bv}{\mathcal U}
\newcommand{\br}{\mathcal R}

\newcommand{\Nn}{\mathbb{N}}
\newcommand{\Rr}{\mathbb{R}}

\newcommand{\Zz}{\mathbb{Z}}
\newcommand{\Cc}{\mathbb{C}}

\newcommand{\si}{\sigma}
\newcommand{\Si}{\Sigma}
\newcommand{\om}{\omega}

\newcommand{\kah}{\text{K\"{a}hler }}
\newcommand{\hor}{H\"ormander }
\newcommand{\dbar}{\bar \partial}
\newcommand{\ddbar}{\partial \bar \partial}

\newcommand{\beq}{\begin{eqnarray*}}
\newcommand{\eeq}{\end{eqnarray*}}
\newcommand{\bpr}{\begin{preuve}}
\newcommand{\epr}{\end{preuve}}

\newenvironment{preuve}[1][]
{\vskip 2mm  {\it \bf Proof#1. }}{$\Box$ \vskip 2mm}

\newtheorem{definition}{Definition}
\newtheorem{theorem}{Theorem}
\newtheorem{lemma}{Lemma}
\newtheorem{proposition}{Proposition}
\newtheorem{corollary}{Corollary}

\newcommand{\fsx}{f_{\sigma, x}}

\newcommand{\gsx}{g_{ x}}
\newcommand{\pgsx}{\partial g_{ x}}
\newcommand{\pfsx}{\partial f_{\sigma, x}}

\newcommand{\rhx}{\Rr H_{x}}
\newcommand{\rhdx}{\Rr H_{2x}}
\newcommand{\rhtx}{\Rr H_{3x}}
\newcommand{\xz}{x}
\newcommand{\sumn}{\sum_{i=0}^{n-1}}
\newcommand{\sumunn}{\sum_{i=1}^{n-1}}

\newcommand{\hld}{H^0(X;L^d)\setminus \widetilde \Delta_d}
\newcommand{\rhld}{\Rr H^0(X;L^d)\setminus  \Rr \widetilde\Delta_d}
\newcommand{\tdelta}{\widetilde \Delta_d}

\newcommand{\zsz}{\Zz/2\Zz}

\newcommand{\rhl}{\Rr H^0(X;L^d)}
\newcommand{\Rhd}{\Rr H^0(X;L^d)}

\newcommand{\unan}{\{1, \cdots, n-1\}}
\newcommand{\cpun}{\Cc P^1}
\newcommand{\rpun}{\Rr P^1}

\title{What is the total Betti number of \\a random real hypersurface?}
\author{Damien Gayet, Jean-Yves Welschinger}

\begin{document}
\large
\maketitle
\centerline{\textbf{Abstract}}
We bound from above  the expected total Betti number of 
a  high degree random real hypersurface  in a smooth real projective manifold. 
This upper bound is deduced from the equirepartition of critical points of 
a real Lefschetz pencil restricted to the complex domain of such a random hypersurface, equirepartition which 
we first establish. Our proofs involve H\"ormander's theory of peak sections as well as the formula of
 Poincar\'e-Martinelli. 


\section*{Introduction}
The topology of real projective manifolds is under study since the nineteenth century, when Axel Harnack
and Felix Klein discovered that the number of  connected components
of the real locus of a smooth real projective curve is bounded from above 
by the sum of its genus and the number of connected components of its complex domain, see 
\cite{Harnack}, \cite{Klein}, while
David Hilbert has devoted his sixteenth problem to such a study. 
Recall that by definition, a real projective manifold $X$ is the vanishing locus in
some complex projective space of a collection of homogeneous polynomials with real coefficients. It
inherits an antiholomorphic involution $c_X$ 
from the ambient complex conjugation. 
The real locus $\Rr X$ is the set of real solutions of the polynomial equations, that 
is the fixed point set of $c_X$. 
Ren\'e Thom  \cite{Thom} later observed  as a consequence
of Smith's theory in equivariant homology, that 
the total Betti number of the real locus of a smooth real projective manifold is actually always
bounded from above by the total Betti number of its complex locus, extending
Harnack-Klein's inequality, see Theorem \ref{theo Smith}. On the other hand,
John Nash proved that every closed smooth manifold can be realized as a component of 
the real locus of a smooth real projective manifold.

Real projective manifolds  achieving the upper bound given by  
Harnack-Klein or Smith-Thom's inequalities are called maximal. Real maximal curves in 
 smooth real projective surfaces appear to be exponentially rare
in their linear system as their degree grows, see \cite{GW}. What is then  the 
expected topology of 
 real hypersurfaces in a given smooth real projective manifold $X$? 
We  tackle  here this question, 
measuring
 the topology of hypersurfaces by the total Betti numbers 
of their real loci.
The answer 
to this question indeed turns out to be only known for  the real projective line 
thanks  to Mark Kac  \cite{Kac}, Michael Shub and Stephen Smale \cite{SS} or Alan Edelman and Eric Kostlan \cite{EK}. 
From these works  follows  that the expected 
number of real roots of a random real polynomial in one variable
and degree $d$ is $\sqrt d$.
We establish here general upper bounds for the expected
total Betti numbers of real hypersurfaces in real projective manifolds. More 
precisely, let $X$ be a 
smooth real projective manifold   of positive dimension $n$
equipped with a real ample line bundle $L$. The growth of the  total Betti number
of complex loci of hypersurfaces  linearly equivalent to $L^d$
is polynomial in  $d$ of  degree $n$, see Lemma \ref{lemma Betti asymptotic}. 
We prove the following, see Theorem \ref{theo principal}
and \ref{theo principal produit}. 
\begin{theorem}\label{theo zero}
Let $(X,c_X)$ be a smooth real projective manifold of dimension $n$ greater than one
 equipped with a Hermitian real line bundle $(L,c_L)$
of positive curvature.  Then, 
the expected total Betti number of real loci of hypersurfaces 
 linearly equivalent to $L^d$ is a $o(d^{n})$.
If $n=2$ or if $X$ is a product
of smooth real projective curves, then 
 it is even a $O(d^{\frac{n}{2}}(\log d)^n)$. 
\end{theorem}
The probability measure that we consider on the complete linear system 
of real divisors associated to $L^d$ is the Fubiny-Study measure arising from the $L^2$-scalar product 
induced by the  Hermitian metric  of positive curvature fixed on $L$, 
 see \S \ref{para 3.1}. When $X$ is one-dimensional, 
upper bounds as the ones given by Theorem \ref{theo zero}
can already be deduced from our work \cite{GW}. 
 In order to prove Theorem \ref{theo zero}, 
we first fix 
a real Lefschetz pencil on $X$, which 
 restricts to a Lefschetz pencil on every generic hypersurface
of $X$. The number of critical points of such a restriction 
has the same asymptotic as the  total Betti number
of the hypersurface, see \S \ref{para asymptotics}.  We then prove that these critical
points get uniformly distributed in $X$ when the degree increases and more precisely that
the expected normalized counting measure supported by these critical points 
converges to the volume form induced by the curvature of the Hermitian bundle $L$,
see Theorem \ref{theo equi real}. The latter  weak convergence 
proved in Theorem \ref{theo equi real} is established outside  of the critical locus of the 
original 
pencil when  $n>2$ and away from  the real locus of $X$. Note 
 that we first prove this equirepartition result over   the complex numbers,
see Theorem \ref{theo equi}. 
In order to deduce Theorem \ref{theo zero} from this equirepartition result,
we observe that 
the total Betti number of real loci of hypersurfaces is  bounded from above
by the number of critical points of the restricted Lefschetz pencil 
in a neighborhood
of the real locus, which we choose of size $\frac{\log d}{\sqrt d}$
thanks to the theory of H\"ormander's peak sections, see \S \ref{para 3}.
The bound $d^{\frac{n}{2}}(\log d )^n$ of Theorem \ref{theo zero}
indeed appears to be the volume of a $\frac{\log d}{\sqrt d}$-neighborhood 
of the real locus for the metric induced by the curvature form of $L^d$.

 Theorems \ref{theo equi} and \ref{theo equi real} on equirepartition of critical points
are  independent of Theorem \ref{theo zero}
which motivated this work.  Note also  that  the expected Euler characteristic of the real locus 
of such random hypersurfaces has been computed in \cite{IP} and \cite{Burgisser}. 
Finally, while we were writing this paper in june 2011, Peter Sarnak informed us that he
is able to prove together with Igor Wigman in a work in progress that the expected number of connected
components of  real curves of degree $ d$ in $\Rr P^2$ 
is even a $O(d)$. 

Our paper is organized as follows. In the first paragraph, we recall few results
about total Betti numbers of real projective manifolds,  critical points of Lefschetz
pencils and their asymptotics. The second paragraph is devoted to Theorem \ref{theo equi}
about  equirepartition of critical points of Lefschetz pencils restricted to complex random 
hypersurfaces. 
The theory of peak sections of \hor  and Poincar\'e-Martinelli's formula play a crucial r\^ole 
in the proof, see \S \ref{para PM formula} and \S \ref{para 3}. Finally in the third paragraph, we first establish the real analogue of Theorem \ref{theo equi},
see Theorem \ref{theo equi real}, and then  deduce  Theorem \ref{theo zero} from it,
that is upper bounds for the expected
 total Betti numbers of the real locus of random real hypersurfaces, see Theorems \ref{theo principal}
and  \ref{theo principal produit}. 

\textit{Aknowledgements.} The research leading to these results has received funding
from the European Community's Seventh Framework Progamme 
([FP7/2007-2013] [FP7/2007-2011]) under
grant agreement $\text{n}\textsuperscript{o}$ [258204], as 
well as from the French Agence nationale 
de la recherche, ANR-08-BLAN-0291-02. 

\tableofcontents

\section{Betti numbers and critical points of Lefschetz pencils}\label{section Betti}
This first paragraph is devoted to Lefschetz pencils, 
total Betti numbers and their asymptotics.
\subsection{Real Lefschetz pencils and Betti numbers}
Let $X$ be a smooth  complex projective  manifold of positive dimension $n$.
\begin{definition}\label{def Lefschetz}
 A Lefschetz pencil on $X$ is a rational map $p: X \dashrightarrow \Cc P^1$ having
only non degenerated critical points and defined by two sections of a holomorphic line bundle
with smooth and transverse vanishing loci.
\end{definition}
We denote by $B$ the base locus of a  Lefschetz pencil $p$ given by Definion \ref{def Lefschetz},
that is the codimension two submanifold of $X$ where $p$ is not defined. 
A Lefschetz pencil without base locus is called a Lefschetz fibration. 
 Blowing up once the base locus of a Lefschetz pencil turns it into a Lefschetz fibration.
When the dimension $n$ of $X$ equals one, 
the base locus is always empty and a Lefschetz fibration is nothing but a branched cover with simple ramifications.
Hence, the following Proposition \ref{prop Euler} extends  to 
Lefschetz fibrations the classical Riemann-Hurwitz formula.
\begin{proposition} \label{prop Euler}
 Let $X$ be a smooth complex projective manifold of positive dimension $n$ equipped with a Lefschetz 
fibration $p: X \rightarrow \Cc P^1$ and let $F$ be a regular fiber of $p$. Then, the Euler characteristics 
of $X$ and $F$ satisfy the relation 
$$ \chi(X) = 2\chi(F) + (-1)^n \#\text{Crit}(p),$$
where $\text{Crit}(p)$ denotes the set of critical points of $p$.
\end{proposition}
\bpr Denote by $\infty = p(F) \in \Cc P^1$ and by $F_0$ the fiber of $p$ associated to a 
regular value $0 \in \cpun \setminus \{\infty\}$. Let $U_0
$ (resp. $U_\infty$) be a neighborhood of $0$ (resp. $\infty$)
in $\cpun$, without any critical value of $p$. 
Since $U_\infty$ (resp. $U_0$) retracts on $F$ (resp. $F_0$),
we know that 
$\chi (p^{-1}(U_0)) = \chi (p^{-1}(U_\infty)) = \chi (F),$
whereas from additivity of the Euler characteristic, $\chi(\overline{ X}) = \chi (X) -2\chi (F)$,
where $\overline{ X}$ denotes the complement $X\setminus p^{-1}(U_0\cup U_\infty)$. 
Without loss of generality, we may assume that in an affine chart $\Cc = \cpun \setminus \{\infty\}$,
$0$ corresponds to the origin, $U_0$ to a ball centered at the origin and $U_\infty$ 
to a ball centered at $\infty$. The manifold $\overline{ X}$ comes then equipped with a function 
$ f: x\in \overline{ X} \mapsto |p(x)|^2 \in \Rr^*_+\subset \Cc$, taking values in a compact
interval $[a,b]$ of $\Rr^*_+$. This function $f$ is Morse and has the same critical points as
$p$, all being of index $n$. Indeed, the differential of $f$ writes 
$ df = p\overline{\partial p}+ \overline{p} \partial p $
and vanishes at $x\in \overline{ X}$ if and only if $\partial p_{|x}$ vanishes. 
Moreover, its second differential is the composition of the differential of the norm $| .|^2 $
with the second differential of $p$. The multiplication by $i$ exchanges  stable and unstable
spaces of these critical points which are non degenerated. Hence, $\overline{ X}$
is equipped with a Morse function $f: X \to [a,b]$ having $\# \text{Crit}(p)$ critical points, 
 all of index $n$. By the Morse Lemma (see \cite{Milnor}),
the topology of $f^{-1} ([a,a+\epsilon]) $ changes, as $\epsilon
$ grows, only at the critical points, where a handle $D^n\times D^n$
of index $n$  is  glued on a submanifold diffeomorphic
to $D^n \times S^{n-1}$. From this Morse theory we deduce that $$\chi(\overline{ X}) = \# \text{Crit}(p) (1- \chi (S^{n-1})) 
= (-1)^n \# \text{Crit}(p)$$ and the result.
\epr    
Recall that a complex projective manifold $X\subset \Cc P^N$ is said to be real when
it is defined over the reals,  as the vanishing locus of a system of polynomial equations
with real coefficients. It inherits then an antiholomorphic involution $c_X:X\to X$, 
which is the restriction of the complex conjugation  $conj: (z_0: \cdots: z_n) \in \Cc P^n \to 
(\overline{ z_0}: \cdots: \overline{ z_n})\in \Cc P^n$. Its fixed point set $\Rr X\subset \Rr P^n$ is called
the real locus of $X$. When $X$ is smooth, the latter is either empty or half-dimensional.
\begin{definition}\label{def real Lefschetz}
 Let $(X,c_X)$ be a  smooth real projective manifold of positive dimension $n$. 
A Lefschetz pencil $p: X\dashrightarrow \Cc P^1$  is said to be real 
iff it satisfies $p\circ c_X = conj \circ p$. 
\end{definition}
Such a real Lefschetz pencil given by Definition \ref{def real Lefschetz}
is then defined by two real sections $\si_0$, $\si_1$ of a holomorphic real line bundle $\pi: (N,c_N) \to (X,c_X)$,
where $\pi \circ c_N = c_X \circ \pi$.  

Now, if $M$ is a smooth manifold of positive dimension $n$, we denote by 
$$b_*(M; \Zz/2\Zz) = \sum^n_{i=0} \dim H_i(M; \Zz/2\Zz)$$ its total Betti number with
  $\Zz/2\Zz$-coefficients.
\begin{lemma}\label{lemma Betti}
Let $M$ be a smooth manifold equipped with a smooth fibration $p: M \to \Rr P^1$
and $F$ be a regular fiber of $p$. Then, the total Betti numbers of  $M$ and $F$ satisfy
$$ b_*(M; \zsz) \leq 4b_*(F; \zsz) + \# \text{Crit}(p).$$
This relation also holds when $M$ is the real
locus of a smooth real projective manifold and $p$
the restriction of a real Lefschetz pencil.
\end{lemma}
\bpr
Denote by $\infty =p(F)\in \rpun$  and by $F_0$ the fiber of $p$ associated to 
a regular value $0\in \Rr P^1\setminus \{\infty\}$. Let $I_0$ (resp. $I_\infty$)
be a neighborhood of $0$ (resp. $\infty$) in $\rpun$, so that 
$I_0$ and $I_\infty$ cover $\rpun$, such that $I_\infty$ contains only regular values of $p$.
Now, set
$U_0 = p^{-1}(I_0)$ and  $U_\infty = p^{-1}(I_\infty)$, so that 
$U_0 \cup U_\infty =M$. 
 From the Mayer-Victoris formula follows that
$$
b_* (M) \leq  b_*(U_0 ) +b_*(U_\infty) + b_*(U_0\cap U_\infty) \leq  b_*(U_0) + 3 b_*( F),
$$
since we may assume  that $U_0\cap U_\infty$ retracts onto two  fibers of $p$.
Now, the restriction of $p$ to $U_0$ is a Morse function  taking values in $I_0$ and
having the same critical points as $p$. 
By the Morse Lemma,  
$b_*(U_0) \leq b_*(F) + \# \text{Crit}(p_{|\Rr X})$.
This proves the first part of Lemma \ref{lemma Betti}. 

If $p: X \dashrightarrow \Cc P^1$ is a real Lefschetz pencil with base locus
$B$, we denote by $ \widetilde X \to X$ the blow-up of $B$ in 
$X$ and by $\widetilde p: \widetilde X \to \cpun$
 the induced Lefschetz fibration. From what has just been proved, we know that
$$ b_*(\Rr \widetilde X; \zsz) \leq 4b_*(\Rr \widetilde F, \zsz) + \# \text{Crit}(\widetilde p_{|\Rr \widetilde X}),$$
where $\widetilde F$ denotes the fiber of $\widetilde p$ associated to $F$.
Moreover the morphism $H_* (\Rr \widetilde X; \zsz) \to H_* (\Rr  X; \zsz) $ 
is onto, since every element of $H_* (\Rr  X; \zsz) $ has a representative transverse to $\Rr B$ and
a proper transform in $\Rr \widetilde X$. It follows that $b_*(\Rr X) \leq b_*(\Rr \widetilde X)$ whereas
the projection $ \Rr \widetilde F \to \Rr F$ is a diffeomorphism. 
\epr
Recall finally the following Theorem \ref{theo Smith} proved by R. Thom in \cite{Thom},
as a consequence of Smith's exact sequence in equivariant homology. 
\begin{theorem}\label{theo Smith}
Let $(X,c_X)$ be a smooth real projective manifold with real locus $\Rr X$. 
Then, the total Betti numbers of $X$ and $\Rr X$ satisfy
$ b_*(\Rr X; \zsz) \leq b_*(X;\zsz)$.  
\end{theorem}
The manifolds for which equality holds in Theorem \ref{theo Smith} are called maximal. For instance,  real projective spaces 
are maximal. When $X$ is one-dimensional and irreducible, Smith-Thom's inequality  given by 
Theorem \ref{theo Smith} reduces to the Harnack-Klein's inequality, up to which 
the number of connected components of $\Rr X$ is bounded from above by $g(X)+1$, where $g(X)$ denotes
the genus of the curve $X$, see \cite{Wilson} and references therein. Real maximal curves in real projective
surfaces turn out to become exponentially rare  in their linear system as their degree grows, see \cite{GW}.

\subsection{Asymptotics}\label{para asymptotics}
Given a holomorphic line bundle $L$ over a smooth complex projective  manifold $X$, we denote, for every
non trivial section $\si$ of $L$, by $C_\si$ its vanishing locus. 
\begin{lemma} \label{lemma Chern}
 Let $L$ be a holomorphic line bundle over a smooth complex projective  manifold $X$ of positive dimension $n$.
For every section $\si$ of $L$ which vanishes transversally, the Chern classes of its vanishing locus $C_\si$ write:
$$ \forall j\in \unan, c_j (C_\si) = \sum ^j _{k=0} (-1)^k c_1(L)^k \wedge c_{j-k} (X)_{|C_\si}\in H^{2j}(C_\si;\Zz).$$
\end{lemma}
In particular, if $(\si_d)_{d>0}$ is a sequence of sections of $L^d$ given by Lemma \ref{lemma Chern},
the Euler characteristic of $C_{\si_d}$ is a polynomial of degree $n$ in $d$ with leading coefficient 
 $(-1)^{n-1} \int_X c_1(L)^n$.
\bpr
The adjunction formula for $X$, $C_\si$ and $L$ writes 
$ c(X)_{|C_\si} = c(C_\si) \wedge c(L)_{|C_\si},$
since the restriction of $L$ to 
$C_\si$ is isomorphic to the normal bundle of $C_\si$ in $X$. 
As a consequence, $c_1(X)_{|C_\si} = c_1(C_\si) + c_1(L)_{|C_\si}$ and
for every $j\in \{2, \cdots, n-1\},$
$$c_j(X)_{|C_\si} = c_j(C_\si) + c_{j-1}(C_\si) \wedge c_1(L)_{|C_\si}.$$ 
Summing up, we get the result.
\epr
\begin{lemma}\label{lemma Betti asymptotic}
 Let $L$ be an ample line bundle over a smooth complex projective  manifold of positive dimension $n$. 
Let $(\si_d)_{d>0}$ be a sequence of sections of $L^d$
vanishing transversally. Then, 
$$ b_*(C_{\si_d};\Zz/2\Zz) = (-1)^{n-1} \chi (C_{\si_d}) + O(1) = \left(\int_X c_1(L)^n\right) d^n
+ O(d^{n-1}).$$ 
\end{lemma}
\bpr
When $d$ is large enough, $L^d$ is very ample and we choose
an embedding of $X$ in $\Cc P^N$, $N>0$, such that $L^d$ coincides with 
the restriction of $\mathcal O _{\Cc P^N}(1)$ to $X$. Then, $C_{\si_d} $
writes $X\cap H$ where $H$ is a hyperplane of $\Cc P^N$. By  Lefschetz's
theorem of  hyperplane sections, for $0\leq i \leq n-1$, 
$\dim H_i(C_{\si_d}; \zsz) = \dim H_i(X;  \zsz)$
and then by Poincar\'e duality, 
$\dim H_ {2n-2-i}(C_{\si_d}; \zsz) = \dim H_{2n-i}(X;  \zsz)$.
Hence 
\beq
b_*(C_{\si_d};\Zz/2\Zz) &=& \dim H_{n-1} (C_{\si_d}; \zsz) + O(1)\\
 &=& (-1)^{n-1} \chi(C_{\si_d}) + O(1).
\eeq
The result now follows  from Lemma \ref{lemma Chern}.
\epr

\begin{proposition}\label{prop critical}
 Let $X$ be a smooth complex projective  manifold of dimension $n$ greater than one
equipped with a Lefschetz pencil $p: X \dashrightarrow \cpun$. Let $L\to X$ 
be a holomorphic line bundle and $\si_d$ be a section of $L^d$ which vanishes 
transversally, where $d>0$. Assume that the restriction of $p$ to $C_{\si_d}$
is Lefschetz. Then, the number of critical points of the restriction $p_{|C_{\si_d}}$
equals $(\int_X c_1(L)^n ) d^n+ O(d^{n-1})$.
\end{proposition}
\bpr
Denote by $\widetilde X$ (resp. $\widetilde C_{\si_d}$) the blow-up of the base locus $B$ (resp. $B\cap C_{\si_d}$)
of $p$ (resp. $p_{|C_{\si_d}}$), so that $\widetilde X$ (resp. $\widetilde C_{\si_d}$) is equipped with a Lefschetz
fibration induced by $p: \widetilde X \to \cpun$ (resp. $p_{|\widetilde C_{\si_d}}:\widetilde C_{\si_d} \to \cpun$). Let
$F$ be a regular fiber of $p$ transverse to $C_{\si_d}$ and $\widetilde F$ be the corresponding fiber in $\widetilde X$. 
By Proposition \ref{prop Euler}, 
$$ (-1)^{n-1} \# \text{Crit}(p_{|C_{\si_d}}) = \chi(\widetilde C_{\si_d}) -2\chi (\widetilde F\cap \widetilde C_{\si_d}).$$
From additivity of the Euler characteristic, we know that
 $\chi(\widetilde C_{\si_d}) = \chi (C_{\si_d}) + \chi (B\cap C_{\si_d})$.
The exceptionnal divisor of $\widetilde C_{\si_d}$ over $B\cap C_{\si_d}$ is indeed a ruled surface
over $B\cap C_{\si_d}$ of Euler characteristic $2\chi(B\cap C_{\si_d})$. 
Likewise, 
$\chi(\widetilde F\cap \widetilde C_{\si_d}) = \chi (F\cap C_{\si_d}),$
since the projection $\tilde F\cap \widetilde C_\si\to F\cap C_\si$
is a diffeomorphism.
The result now follows from  Lemma \ref{lemma Chern}, which  provides the equivalents
$$ \chi(C_{\si_d}) \sim_{d\to \infty} (-1)^{n-1}\left(\int_X c_1(L)^n\right)d^n,$$
$$ \chi(B\cap C_{\si_d}) \sim_{d\to \infty} (-1)^{n-3}\left(\int_B c_1(L)^{n-2}_{|B}\right)d^{n-2} \text{ and}$$
$$ \chi(F\cap C_{\si_d}) \sim_{d\to \infty} (-1)^{n-2}\left(\int_F c_1(L)^{n-1}_{|F}\right)d^{n-1}.$$ 
\epr

\section{Random divisors and distribution of critical points}\label{section random divisors}
Let $X$ be a smooth complex projective  manifold equipped with a Lefschetz pencil. 
The restriction of this pencil to a generic smooth hypersurface 
$C$ of $X$ is a Lefschetz pencil of $C$. The aim of this paragraph is to prove
the equidistribution in average of critical points of such a restriction to a random hypersurface $C$ 
of large degree, see Theorem 
\ref{theo equi}. The estimations of the total Betti number of real hypersurfaces
will be obtained in paragraph \ref{section real} as a consequence of a real
analogue of this Theorem \ref{theo equi}, see Theorem \ref{theo equi real}. 

We first formulate this equirepartition Theorem \ref{theo equi}, then introduce the main
ingredients of the proof, namely  Poincar\'e-Martinelli's formula and 
 H\"ormander's peak sections. Finally, we prove Theorem \ref{theo equi}.
Note that this paragraph is independent of the remaining part of the paper, 
it does not involve any real geometry. 
\subsection{Notations and result}\label{para notations}
Let $X$ be a smooth complex projective  manifold of positive dimension $n$ equipped with 
a Lefschetz pencil $ p: X \dashrightarrow \Cc P^1$ with base locus $B\subset X$. 
Let $L\to X$ be a holomorphic line bundle equipped with a Hermitian metric $h$ of
positive curvature $\omega \in \Omega^{(1,1)}(X;\Rr)$. The latter is defined
in the neighborhood of every point $x\in X$ by the relation 
$ \omega = \frac{1}{2i\pi} \partial \dbar \log h(e,e)$, where $e$ is
a local non vanishing holomorphic section of $L$ defined in the neighborhood of $x$.
The curvature form induces a K\"ahler metric on $X$ and we denote by $dx = \frac{\omega^n}{\int_X \omega^n}$
its associated normalized volume form. For every integer $d>0$, we denote by 
$h^d$ the induced Hermitian metric on the bundle $L^d$ and by $\langle\  \rangle$ the induced  $L^2$-Hermitian  product
 on the space $H^0(X;L^d)$ of global sections of $L^d$. This product is defined by the relation
 $$ (\sigma, \tau) \in H^0(X;L^d)\times H^0(X;L^d) \mapsto \langle \sigma, \tau \rangle = \int_X h^d(\sigma, \tau) dx \in \Cc.$$
 Denote by $N_d$ the dimension of $H^0(X;L^d)$ and by $\mu$ its Gaussian measure, defined by the relation 
$$ \forall A \subset H^0(X;L^d), \mu(A) = \frac{1}{\pi^{N_d}}\int_A e^{-||\sigma||^2}d\sigma,$$
where  $||\sigma||^2 = \langle \sigma , \sigma \rangle$ and $d \sigma$ denotes  the Lebesgue measure associated
to $\langle \ \rangle$. 
Denote by $\Delta_d\subset H^0(X;L^d)$ the discriminant locus, that is  the set of  sections of $H^0(X;L^d)$
which do not vanish tranversally.
Likewise, denote by $\widetilde \Delta_d \in H^0(X;L^d)$ the union of $ \Delta_d$ with the set of 
sections $\sigma\in H^0(X;L^d)$ such that either the restriction of 
 $p$ to  $C_\sigma$  is not Lefschetz, or this vanishing locus $C_\si$ 
meets the critical set $\text{Crit}(p)$. 
By Bertini's theorem (see for example Theorem 8.18 of \cite{Hartshorne}), $\widetilde \Delta_d $ is 
a hypersurface of $H^0(X;L^d)$ as soon as $d$ is large enough, which will be assumed throughout this article.

 For every section $\sigma \in H^0(X;L^d)\setminus \widetilde \Delta_d$, denote by $\mathcal R_\si$ the set of 
 critical points of the restriction $p_{|C_\sigma}$ of $p$ to $C_\sigma$, so that by Proposition \ref{prop critical}, 
 the cardinal $\# \mathcal R_\sigma$ of this set is equivalent to $(\int_X \omega^n) d^n$ as $d$ grows
 to infinity. For every $x\in X$, we finally denote by $\delta_x
$ the Dirac measure $\chi\in C^0(X,\Rr) \mapsto \chi(x)\in \Rr$.
\begin{definition}\label{def Dirac} For every $\sigma \in H^0(X;L^d)\setminus \widetilde \Delta_d$, 
the measure $\nu_\sigma = \frac{1}{\#\mathcal R_d}\sum _{x\in \mathcal R_\sigma} \delta_x$ is called the probability measure 
of $X$ carried by the critical points of $p_{|C_\sigma}$. 
\end{definition}
Our goal in this paragraph is to prove the following Theorem \ref{theo equi} which asymptotically computes 
the expected  probability measure given by Definition \ref{def Dirac}.
\begin{theorem}\label{theo equi}
Let $X$ be a smooth complex projective  manifold of dimension $n$ greater than one equipped with a 
 Lefschetz pencil $ p: X \dashrightarrow \Cc P^1$ with critical locus $\text{Crit}(p)$. 
Let $L\to X$ be a holomorphic line bundle equipped with a Hermitian metric $h$ of
positive curvature $\omega$. Then, for every function $\chi: X \to \Rr$ of class $C^2$ such that,
 when $n>2$, the support of $\partial \dbar \chi$ 
is disjoint  from $\text{Crit}(p)$, we have
$ \lim_{d\to \infty} E(\langle \nu_\sigma, \chi \rangle) = \int_X \chi dx,$
where $E(\langle \nu_\sigma, \chi \rangle) = 
\int_{H^0(X;L^d)\setminus \widetilde \Delta_d} \langle \nu_\sigma, \chi \rangle d\mu(\sigma)$.
\end{theorem}
Note that similar results as Theorem \ref{theo equi}, on equirepartition of critical points of sections,
have been obtained in \cite{DSZ1}, \cite{DSZ2}
by M. Douglas, B. Shiffman and S. Zelditch.

Note also that the equirepartition Theorem \ref{theo equi}, as well as Theorem  \ref{theo equi real},
is local in nature and does not depend that much on a Lefschetz pencil.
Any  local holomorphic Morse function could be used instead of a Lefschetz pencil,
leading to the same proof and conclusions.


\subsection{Poincar\'e-Martinelli's formula and adapted atlas}\label{Adapted atlas}
\subsubsection{Adapted atlas and  associated relative trivializations}\label{para atlas}

\begin{definition} \label {def adapted}Let $X$ be a smooth complex projective manifold of positive dimension $n$ 
equipped with a Lefschetz pencil $ p: X \dashrightarrow \Cc P^1$. An atlas $\mathcal U$ of $X$ 
  is said to be \textit{adapted} to $p$ iff for every open set $U\in \mathcal U$, the restriction of $p$
  to $U$ is conjugated to one of the following three models in the neighborhood of the origin
  in $\Cc^n$:
  \begin{itemize}
  \item [(r)] $(z_1, \cdots , z_n) \in \Cc^n \mapsto z_n \in \Cc$
  \item [(b)] $(z_1, \cdots , z_n) \in \Cc^n \setminus \Cc^{n-2}  \mapsto [z_{n-1}: z_n] \in \Cc P^1$
  \item [(c)] $(z_1, \cdots , z_n) \in \Cc^n  \mapsto z_1^2 + \cdots +z^2_n \in \Cc$
  \end{itemize} 
  \end{definition}
  Every atlas of $X$ becomes adapted in the sense of Definition \ref{def adapted}
after refinement. Let $x$ be a point in  $X$.
  If $x$ is a regular point of $p$, by the implicit function theorem it has a neighborhood biholomorphic to the model
  $(r)$ of Definition \ref{def adapted}. If $x$ is a base point 
(resp. a critical point), it has by definition (resp. by the holomorphic
  Morse Lemma) a neighborhood biholomorphic to the model $(b)$ (resp. $(c)$) of Definition \ref{def adapted}.
 
 In the model $(r)$, the vertical tangent bundle $\ker (dp)$ is trivialized by the vector fields 
 $\frac{\partial}{\partial z_1}, \cdots, \frac{\partial}{\partial z_{n-1}} $ of $\Cc^n$. In the model $(b)$, it is trivialized outside of
 the base locus by the vector fields 
 $z_{n-1} \frac{\partial}{\partial z_{n-1}}+ z_{n}\frac{\partial}{\partial z_{n}},  \frac{\partial}{\partial z_1}, \cdots, \frac{\partial}{\partial z_{n-2} }$ of $\Cc^n$.
 In the model $(c)$, when $n= 2$, it is trivialized outside of the critical point by the vector field 
 $ z_1\frac{\partial}{\partial z_2}-z_2\frac{\partial}{\partial z_{1}} $ of $\Cc^2$.
\begin{definition} \label {def trivialization}  
Let  $X$ be a smooth complex projective manifold of positive dimension $n$ 
equipped with a Lefschetz pencil $ p: X \dashrightarrow \Cc P^1$ and an adapted atlas $\mathcal U$. A relative trivialization 
associated to $\mathcal U$ is the data, for every open set $U\in \mathcal U$, of $n-1$ vector fields on $U$ corresponding
to the vector fields 
$\frac{\partial}{\partial z_1}, \cdots, \frac{\partial}{\partial z_{n-1}} $ in the model $(r)$, to 
 $z_{n-1}\frac{\partial}{\partial z_{n-1}}+ z_{n}\frac{\partial}{\partial z_{n}}, \frac{\partial}{\partial z_1}, \cdots, \frac{\partial}{\partial z_{n-2} }$ in the model
  $(b)$ and to 
  $ z_1\frac{\partial}{\partial z_2}-z_2\frac{\partial}{\partial z_{1}} $ in the model $(c)$ when $n=2$.
 \end{definition} 
Note that in the model $(c)$ given by Definition \ref{def adapted} the vertical tangent bundle $\ker (dp)$
restricted to $\Cc^n\setminus \{0\}$ is isomorphic to the pullback of the cotangent bundle of $\Cc P^{n-1}$ by the projection
$\pi: \Cc^n\setminus \{0\} \to \Cc P^{n-1}$. Indeed, the fibers of this 
vertical tangent bundle are the kernels of the 1-form $\alpha = \sum_{i=1}^{n}z_i dz_i$,
so that the restriction map induces  an isomorphism $(\Cc^n)^*/<\alpha> \cong (\ker dp)^*$. 
But the canonical identification between $\Cc^n$ and $(\Cc^n)^*$ gives an isomorphism between the bundle $(\Cc^n)^*/<\alpha>$ and 
$\pi^*(T P ^{n-1})$ over $\Cc^n \setminus \{0\}$. By duality, we get the isomorphism $\ker (dp) \cong \pi ^*T^*\Cc P^{n-1}$. 
When $n>2$, we no more see trivialisations
of this bundle over $\Cc^n \setminus \{0\}$ and thus restrict ourselves to $n=2$
for the model $(c)$ in Definition \ref{def asso}.

\begin{definition}\label{def asso}
Let $X$ be a smooth complex projective  manifold of positive dimension $n$ equipped with 
a Lefschetz pencil $ p: X \dashrightarrow \Cc P^1$ and with 
  a holomorphic line bundle $L\to X$. An atlas $\mathcal U$ is said to be adapted to $(p,L)$
 if it is adapted to $p$ in the sense of Definition \ref{def adapted} and if for every open set $U\in \bv$, 
 the restriction of $L$ to $U$ is trivializable. A relative trivialization associated to $\bv$
 is a relative trivialization in the sense of Definition \ref{def trivialization} together with 
 a trivialization $e$ of $L_{|U}$,  for every open 
 set $U$ of $\mathcal U$.
 \end{definition}
\subsubsection{ Poincar\'e-Martinelli's formula}\label{para PM formula}
Let $X$ be a smooth complex projective  manifold of positive dimension $n$ equipped with 
a Lefschetz pencil $ p: X \dashrightarrow \Cc P^1$ of critical locus $\text{Crit}(p)$. Let $L\to X$ 
be an ample holomorphic line bundle equipped with a Hermitian metric $h$ of positive curvature $\omega\in \Omega^{(1,1)}(X;\Rr)$.
Let $\mathcal U$ be an atlas of $X$ adapted to $(p,L)$ and $(v_1, \cdots, v_{n-1}, e)$ be an associated relative trivialization 
given by Definition \ref{def asso}. Let $U$ be an element of $\mathcal U$. For every section $\sigma \in H^0(X;L^d)$, 
we denote by $f_{\sigma, U}: U \to \Cc$ the holomorphic function defined by the relation $\si_{|U} = f_{\si,U} e^d_{|U}$. 
When $\si \notin \widetilde \Delta_d$, the set $\br_\si\cap U$ coincides by definition with the
transverse intersection of the hypersurfaces $\{f_{\si,U} =0\}, \{\partial f_{\si,U}(v_1) =0\}, 
\cdots,\{\partial f_{\si,U}(v_{n-1}) =0\}$. For every function $\chi: X \to \Rr$ with compact support in $U$, 
 Poincar\'e-Martinelli's formula (see \cite{Griffiths-King}) then writes here:
\begin{equation}\label{PM} \langle \nu_\si, \chi \rangle = (\frac{i}{2\pi})^n \frac{1}{\# \br_\si}
\int_X \lambda_U \ddbar  \chi\wedge
(\ddbar \lambda_U)^{n-1},\tag{PM}
\end{equation} 
where $$\lambda_U = \log \left(d^2 |f_{\si,U}|^2 + \sum^{n-1}_{i=1} |\partial f_{\si,U}(v_i)|^2\right).$$
Note that in this definition of $\lambda_U$, we have chosen for convenience to use  $df_{\si,U}$
instead of $f_{\si,U}$ as the first function.
This formula of Poincar\'e-Martinelli computes the integral of $\chi$ for the measure $\nu_\si$ introduced in Definition \ref{def Dirac};
its left hand side does not involve any trivialization of $L$ over $U$, 
contrary to the right hand side. It makes it possible to  estimate 
the expectation of the random variable $\langle \nu_\si,\chi\rangle$. 
However, it appears to be useful for this purpose to choose an appropriate trivialization of $L$
in the neighborhood of every point $x\in X$,  whose norm reaches
a local maximum at $x$ where it equals one. We are going to make such a choice instead of 
the trivialization $e$ defined on the whole $U$, as discussed in the
following Proposition \ref{prop PMx} and \S \ref{para 3}.
\begin{proposition} \label{prop PMx}
Let $X$ be a smooth complex projective  manifold of positive dimension $n$ equipped with a 
 Lefschetz pencil $ p : X \dashrightarrow \Cc P^1$. Let $L$ 
be an ample holomorphic line bundle equipped with a Hermitian metric $h$ of positive curvature $\omega$.
 Let $U$ be an element of an atlas adapted to $(p,L)$ and $(v_1, \cdots , v_{n-1},e)$ be an associated
 relative trivialization. Finally, let $(g_x)_{x\in U}$ be a family of germs of holomorphic functions
 such that $g_x$ is defined in a neighborhood of $x$ and
 $\Re g_x(x) = -\log h^d(e^d,e^d)_{|x} $. 
Then, for every function $\chi : X \to \Rr$ of class $C^2$ with support 
in $U$, disjoint from the critical set of $p$ when $n>2$, and for
every $\si \in H^0(X;L^d) \setminus \widetilde \Delta_d$, we have:
\beq
 \langle \nu_\si, \chi \rangle = \frac{1}{\# \br_\si}d^n
\int_X  \chi  \om^n + \frac{1}{\# \br_\si} \sum_{k=0}^{n-1} (\frac{i}{2\pi})^{n-k}d^k
\int_X  \ddbar \chi \wedge \omega^k\wedge \lambda_x (\ddbar \lambda_x)^{n-1-k},
\eeq
 where $$\lambda_x =  
 \log \left(d^2 |f_{\si,x}|^2 + \sum^{n-1}_{i=1} |\partial f_{\si,x}(v_i)+f_{\si,x} \partial g_x(v_i)|^2\right)$$
 and $\si = f_{\si,x} \exp(g_x) e^d$ in the neighborhood of every point $x\in U$.
 \end{proposition}
The condition on $(g_x)_{x\in U}$ in Proposition \ref{prop PMx} ensures that
$\exp(g_x) e^d$ is a holomorphic trivialization of norm one at $x$, so that
$h^d(\si,\si)_{|x}$ coincides with $|f_{\si,x}(x)|^2$. 
The point $x$ in $\lambda_x$ is a parameter and not a variable, so that $\lambda_x$ reads
in the neighborhood of $x$  as a function 
$$z \mapsto   
\log \left(d^2 |f_{\si,x}|^2(z) + \sum^{n-1}_{i=1} |\partial f_{\si,x}(v_i)_{|z}+f_{\si,x}(z) \partial g_x(v_i)_{|z}|^2\right)$$
and $\ddbar \lambda_x$ in the formula given by Proposition \ref{prop PMx}
stands for its second derivative computed at the point $x$.
Note that if $\mathcal U $ is a locally finite atlas adapted to $(p,L)$,  and if $(\rho_U)_{U\in \mathcal U}$
is an associated  partition of unity, then for every function $\chi: X \to \Rr$ of class $C^2$, 
with support disjoint from the critical locus of $p$ when $n>2$, and for every open set $U\in \mathcal U$, 
the function $\chi_U = \rho_U \chi$ satisfies the hypotheses of Proposition \ref{prop PMx},
while $\chi = \sum_{U\in \mathcal U} \chi_U$ and  $ \langle \nu_\si, \chi\rangle  =  
\sum_{U\in \mathcal U} \langle \nu_\si, \chi_U\rangle$.
\bpr Let $\si \in H^0(X;L^d)\setminus \widetilde \Delta_d $ and $x\in U\setminus \text{Crit}(p)$. 
By definition, $f_{\si,U} = f_{\si,x} \exp (g_{x})$ and for every 
$1\leq i \leq n-1$, 
$$ \partial f_{\si,U}(v_i)= (\partial f_{\sigma, x}(v_i) + f_{\sigma, x} \pgsx(v_i)) \exp(g_{x}).$$
As a consequence, 
$\lambda_U  = \Re \gsx+ \lambda_{x}$, so that at the point $x$, $ \lambda_U(x) =-\log h^d(e^d,e^d)(x) +  
\lambda_{x}(x)$.
Since $\ddbar \Re \gsx$  vanishes, the equality $\ddbar \lambda_U = \ddbar  \lambda_{x}$
holds  in a neighborhood of $x$, 
Hence,  formula (\ref{PM}) rewrites 
\beq
\langle \nu_\si, \chi \rangle &=& \frac{i^n}{(2\pi)^n \# \br_\si}
\int_X \ddbar \chi  
\left (-\log h^d(e^d,e^d) + \lambda_{x _0}\right)   (\ddbar \lambda_{U})^{n-1}
\\
&=& \frac{i^{n-1}d}{(2\pi)^{n-1} \# \br_\si} \int_X \ddbar \chi 
\wedge  \omega\wedge 
\lambda_U  (\ddbar \lambda_{U})^{n-2}\\
&&+ \frac{i^n}{(2\pi)^n \# \br_\si}\int_X  \ddbar \chi \wedge
 \lambda_{x _0}  (\ddbar \lambda_{x})^{n-1}
\eeq
The first part of the latter right hand side follows from the relation 
$\ddbar (\lambda_U (\ddbar \lambda_{U})^{n-2}) = (\ddbar \lambda_{U})^{n-1}$, 
the curvature equation $\omega = \frac{i}{2\pi d}\ddbar (-\log h^d(e,e)) $ and Stokes's theorem.
The second part of this right hand side comes from 
 $\ddbar \lambda_U = \ddbar  \lambda_{x}$. Applying this procedure $(n-1)$ times, we deduce by induction
and Stokes's theorem the relation 
\beq
 \langle \nu_\si, \chi \rangle &=& \frac{1}{\# \br_\si}d^n
\int_X  \chi  \om^n\\ &&+ 
\frac{1}{\# \br_\si} \sum_{k=0}^{n-1} (\frac{i}{2\pi})^{n-k}d^k
\int_X  \ddbar \chi \wedge \omega^k\wedge
\lambda_{x} (\ddbar \lambda_{x})^{n-1-k}.
\eeq
\epr
\begin{corollary}\label{cor expectation}
Under the hypotheses of Proposition \ref{prop PMx}, 
\beq
 E(\langle \nu_\si, \chi \rangle) &= &\frac{1}{\# \br_\si}d^n
\int_X  \chi  \om^n  \\ && +\frac{1}{\# \br_\si} \sum_{k=0}^{n-1} (\frac{i}{2\pi})^{n-k}d^k 
\int_X  \ddbar \chi \wedge \omega^k\wedge 
\int_{H^0(X;L^d)\setminus \widetilde \Delta_d}\lambda_x(\ddbar \lambda_x)^{n-1-k} d\mu(\si).
\eeq
\end{corollary}
\bpr 
The result follows by integration over 
$H^0(X;L^d)\setminus \widetilde \Delta_d$
of the relation  given by Proposition \ref{prop PMx}.
\epr
\subsection{H\"ormander's peak sections }\label{para 3}
Let $L$ be a holomorphic line bundle over a smooth complex projective  manifold, equipped with
a Hermitian metric $h$  of positive curvature $\omega $. Let $x$ be a point of $X$. 
There exists, in the neighborhood of $x$, a holomorphic trivialization $e$ of $L$
such that the associated potential $\phi =-\log h(e,e)$ reaches a local minimum at $x$ with
Hessian  of type $(1,1)$. The latter coincides, by definition, with $\omega(.,i.)$. The H\"ormander
$L^2$-estimates makes it possible, for all $d>0$ and maybe after modifying a bit $e^d$ in $L^2$-norm, 
to extend $e^d$ to a global section $\si$ of $L^d$. The latter is called peak section of
H\"ormander, see Definition \ref{def peak}. Moreover, 
G. Tian (Lemma 1.2 in \cite{Tian}) showed that this procedure can be applied to produce 
global sections whose Taylor expansion at $x$ can be controlled at every order, as long
as $d$ is large enough. We recall this result in the following Lemma \ref{lemma peak} where 
we denote, for every $r>0$, by $B(x,r)$ the ball centered at $x$ of radius $r$ in $X$.
\begin{lemma}\label{lemma peak}[See \cite{Tian}, Lemma 1.2]
Let $(L,h)$ be a holomorphic Hermitian line bundle of positive curvature $\omega$ over
a smooth complex projective manifold $X$. Let $x\in X$, $(p_1, \cdots, p_n)\in \Nn^n$
and $p' >p_1+\cdots + p_n$. There exists $d_0\in \Nn$ such that for every $ d>d_0$, the bundle
$L^d$ has a global holomorphic section $\si$ satisfying $\int_X h^d(\si,\si) dx = 1$
and $$ \int_{X\setminus B(x,\frac{\log d }{\sqrt d})} h^d(\si,\si) dx = O(\frac{1}{d^{2p'}}).$$
Moreover, if $z= (z_1, \cdots, z_n)$ are  local coordinates in the neighborhood 
of $x$, we can assume that in a neighborhood of $x$,
$$ \si(z) = \lambda (z_1^{p_1}\cdots z_n^{p_n} + O(|z|^{2p'}))e^d (1+ O(\frac{1}{d^{2p'}} )), $$
where 
$$ \lambda^{-2}= \int_{B(x,\frac{\log d }{\sqrt d})}| z_1^{p_1}\cdots z_n^{p_n}|^2
h^d(e^d,e^d) dx$$ and $e$ is a trivialization of $L$ in the neighborhood of $x$
whose potential $\phi = -\log h (e,e)$ reaches a local minimum at $x$ with
Hessian $\omega(.,i.)$.
\end{lemma} 
\begin{definition}\label{def peak}
We call H\"ormander's peak section of the ample line bundle $L^d$ over the smooth complex
projective manifold $X$ any section given by Lemma \ref{lemma peak} with
$p_1 = \cdots = p_n = 0$ and $p'>1$, where $d$ is large enough.
\end{definition}
Note that such a peak section $\si_0$ given by Definition \ref{def peak}
has its norm concentrated in the neighborhood of the point $x$ given
by Lemma \ref{lemma peak}, so that it is close to the zero section outside of a
 $\frac{\log d }{\sqrt d}$-ball. Moreover, the derivatives and second derivatives
of $\si_0$ at $x$ 
vanish and the value of $\lambda$ at $x$ is equivalent to $\sqrt{(\int_Xc_1(L)^n)d^{n}}$
as $d$ grows to infinity, see Lemma 2.1 of \cite{Tian}.

Note also that if the coordinates $(z_1, \cdots, z_n)$ in Lemma \ref{lemma peak} are
orthonormal at the  point $x$, then two sections given by this lemma 
for different values of $(p_1, \cdots, p_n)$ are asymptotically orthogonal, 
see Lemma 3.1 of \cite{Tian}.
\subsubsection{Evaluation of the two-jets of sections}\label{II 3.2 jets}
Again, let $(L,h)$ be a holomorphic Hermitian line bundle of positive curvature 
$\omega$  over a smooth  $n$-dimensional complex  projective manifold $X$.
Let $x$ be a point of $X$ and $d>0$. We denote by $H_{x}$
the kernel of the evaluation map 
$ \si \in H^0(X;L^d) \mapsto \si(x) \in L^d_{x}$, where $L^d_{x}$
denotes the fiber of $L^d$ over the point $x$. 
Likewise, we denote by $H_{2x}$ the kernel  of the map 
$ \si \in H_{x} \mapsto \nabla\si_{|x} \in T^*_{x}X \otimes L^d_{x}$.
This map does not depend on a chosen connection $\nabla$ on $L$. 
Denote by $H_{3x}$ the kernel  of the map 
$ \si \in H_{2x} \mapsto \nabla^2\si_{|x} \in Sym^2 (T^*_{x}X) \otimes L^d_{x}$.
We deduce from these the jet maps:
\beq
eval_{x}: \si \in H^0(X;L^d)/H_{x} & \mapsto &\si(x) \in L^d_{x},\\
eval_{2x}: \si \in H_{x}/H_{2x} &\mapsto &\nabla \si_{|x} \in T^*_{x}X\otimes L^d_{x},\\
\text{and }  
eval_{3x}: \si \in H_{2x}/H_{3x} & \mapsto &
\nabla^2\si_{|x} \in Sym^2 (T^*_{x}X)\otimes L^d_{x}.
\eeq
When $d$ is large enough, these maps are isomorphisms between 
finite dimensional normed vector spaces. We estimate the norm of these
isomorphisms in the following Proposition \ref{prop norm}, closely following
\cite{Tian}.
\begin{proposition}\label{prop norm}
Let $L$ be a holomorphic Hermitian line bundle of positive curvature over a smooth $n$-dimensional
 complex  projective
manifold $X$. Let $x$ be a point of $X$. Then the maps
$d^{-\frac{n}{2}} eval_{x}$, $d^{-\frac{n+1}{2}} eval_{2x}$ and $ d^{-\frac{n+2}{2}} eval_{3x}$ as
well as their inverse have norms and determinants bounded from above independently of $d$ as
long as $d$ is large enough.
\end{proposition}
Note that Proposition \ref{prop norm} provides an asymptotic result while the condition that $d$ 
be large  ensures  that the three maps are invertible.
\bpr Let $\si_0$ be a peak section of \hor given by Definition \ref{def peak}. By Lemma 2.1 of
\cite{Tian}, $d^{-n}h^d(\si_0, \si_0)_{|x}$ converges to a positive constant as
$d$ grows to infinity. Let $\si_0^{H_{x}}$ be the orthogonal projection of $\si_0$ onto
$H_{x}$. The Taylor expansion of $\si_0^{H_{x}}$ does not contain any constant
term, so that by Lemma 3.1 of \cite{Tian} (see also Lemma 3.2 in \cite{Ruan}), 
the Hermitian product $ \langle \si_0, \frac{\si_0^{H_{x}}}{||\si_0^{H_{x}}||} \rangle$ is a $ O(\frac{1}{d})$, where
$||\si_0^{H_{x}}||^2 = \langle \si_0^{H_{x}} , \si_0^{H_{x}} \rangle$ denotes the $L^2$-norm of $\si_0^{H_{x}}$.
 From the vanishing of
the  product $\langle \si_0 - \si_0^{H_{x}}, \si_0^{H_{x}} \rangle$ we deduce
that $ ||\si_0^{H_{x}}||$ is a $O(\frac{1}{d})$. It follows that the norm of $\si_0 - \si_0^{H_{x}}$ 
equals $1+O(\frac{1}{d^2})$ and we set 
$$\si_0^\perp = \frac{\si_0 - \si_0^{H_{x}}}{||\si_0 - \si_0^{H_{x}}||}.$$
As a consequence, $d^{-n} h^d (\si_0^\perp, \si_0^\perp )_{|x}$ converges to a positive
constant as $d$ grows to infinity. Hence, $d^{-\frac{n}{2}} eval_{x}$ as well as its inverse,
has  norm and determinant bounded when $d$ is large enough. 
The two remaining assertions of  Proposition \ref{prop norm} follow along the same lines. 
For $i \in \{ 1, \cdots n  \}$, let $\si_i$ be a section given by Lemma \ref{lemma peak}, with $p'=2$, 
$p_i = 1$ and $p_j = 0$ if $j\not= i$, $i \in \{ 1, \cdots n  \}$, and where the local coordinates
$(z_1, \cdots , z_n) $ are orthonormal at the point $x$. By Lemma 2.1 of \cite{Tian}, 
for $i \in \{ 1, \cdots n  \}$, $d^{-(n+1)} h^d(\nabla \si_i, \nabla \si_i)_{|x}$ converges
to a positive constant as $d$ grows to infinity. The sections $\si_i$, $i \in \{ 1, \cdots n  \}$,
belong by construction to $H_{x}$ and we set as before 
$$\si_i^\perp = \frac{\si_i - \si_i^{H_{2x}}}{||\si_i - \si_i^{H_{2x}}||},$$
where $\si_i^{H_{2x}}$ denotes the orthogonal projection of $\si_i$ onto $H_{2x}$. 
We deduce  as before from Lemma 3.1 of \cite{Tian} that 
$d^{-(n+1)} h^d(\nabla \si_i^\perp, \nabla \si^\perp_i)_{|x}$ converges
to a positive constant when $d$ grows to infinity and that $h^d(\nabla \si_i^\perp, \nabla \si^\perp_j)_{|x}=0$
if $i\not=j$. The norms of the sections $\si_i^\perp $, $i \in \{ 1, \cdots n  \}$, all equal one but these sections
are not a priori
 orthogonal. However, by Lemma 3.1 of \cite{Tian}, the products 
$\langle \nabla \si_i^\perp, \nabla \si_i^\perp \rangle$ are $O(\frac{1}{d})$ if $j\not=i$, so that
asymptotically, the basis is orthonormal. We deduce that $d^{-\frac{n+1}{2}} eval_{2x}$ 
and its inverse are of norms and determinants bounded as $d$ is large enough. 
The last case follows along the same lines.
\epr

\subsection{Proof of Theorem \ref{theo equi}}\label{subsection proof}
\begin{proposition}\label{prop bounded} Under the hypotheses and notations of Proposition \ref{prop PMx}, 
for all $k \in \{0, \cdots , n-1\}$ and $x\in U$, the integral 
$$\frac{||v_1||^2}{d^k} \int_{H^0(X;L^d)\setminus \widetilde \Delta_d}
|| \lambda_x (\ddbar \lambda_x)^{k} || d\mu(\si)$$ is uniformly bounded by a $O((\log d)^2)$
on every compact subset of $X\setminus \text{Crit}(p)$ when $n>2$ and on the whole $X$  when $n\leq 2$.
\end{proposition}
The norms $|| \ ||$ appearing in the statement of Proposition \ref{prop bounded} are induced  
by the \kah metric of $X$ on elements and 
 $2k$-linear forms of $T_xX$, where $x \in X$. 
 It follows from Definition \ref{def adapted} that $||v_1||$ may vanish in the models 
 $(b)$ and $(c)$. 
 Before proving Proposition \ref{prop bounded}, for which we will spend the whole paragraph, let us first 
deduce a proof of Theorem \ref{theo equi}.
\bpr[ of Theorem \ref{theo equi}]Let $\chi: X \to \Rr$ be a function of class $C^2$ such that the support $K$ of
$\ddbar \chi$ be disjoint from $\text{Crit}(p)$ when $n>2$. Choose a finite atlas $\mathcal U$ adapted 
to $(p,L)$ given by Definition \ref{def asso}, such that when $n>2$, $K$ be covered in 
$X\setminus \text{Crit}(p)$ by elements of $\mathcal U$. Let $U$ be such an element of $\mathcal U$
and $(v_1, \cdots, v_n, e)$ be an associated relative trivialization given by  Definition \ref{def asso}.
Without loss of generality, we can assume that $\chi$ has  support in $U$. 
By Proposition \ref{prop PMx} and with the notations introduced there, 
when $d$ is large enough, the expectation $E(\langle \nu_\si, \chi\rangle)$
equals 
$$ \frac{1}{\# \br_\si}d^n
\int_X  \chi  \om^n + \frac{1}{\# \br_\si} \sum_{k=0}^{n-1} (\frac{i}{2\pi})^{n-k}d^k
\int_X  \ddbar \chi \wedge \omega^k\wedge \int_{H^0(X;L^d)\setminus \widetilde \Delta_d} \lambda_x (\ddbar \lambda_x)^{n-1-k} d\mu(\si).$$
By Proposition \ref{prop critical}, $\# \br_\si $ is equivalent to $(\int_X \omega^n) d^n$ as 
$d$ grows to infinity, so that the first term converges to $\int_X \chi dx$. By Proposition \ref{prop bounded},
the last integral over $H^0(X;L^d)\setminus \widetilde \Delta_d$ is a $O(\frac{(\log d)^2}{d||v_1||^2})$, since we integrate on the support of $\ddbar \chi$ which is disjoint
of $\text{Crit}(p)$ when $n>2$. The result now follows from the fact that the function $||v_1||^{-2}$ is integrable
over $X$. 
\epr

\subsubsection{Proof of Proposition \ref{prop bounded} outside of the base and critical loci of $p$}\label{para 4.1}
Recall that $X$ is equipped with an atlas $\mathcal U$ adapted to $(p,L)$ and with an associated relative trivialization.
The compact $K$ given by Proposition \ref{prop bounded} is covered by 
a finite number of elements of $\mathcal U$. Moreover, we can assume that these elements are all disjoint
from the critical set $\text{Crit}(p)$ when $n>2$. Let $U$ be such an element ; it is either of type $(r)$
given by Definition \ref{def adapted}, or of type $(b)$ or $(c)$. Let us prove now Proposition \ref{prop bounded}
in the case  $U$ be of type $(r)$ and postpone the remaining cases to \S \ref{para b}.

Let $x$ be a point of $K\cap U$. For every $\si \in H^0(X;L^d)\setminus \widetilde \Delta_d$, 
define $h_0 = df_{\sigma, x}$ and for $i\in \unan$, $h_i = \partial f_{\sigma, x}(v_i) + df_{\sigma, x} \pgsx(v_i)$, 
so that 
$$\partial h_i = \partial (\pfsx(v_i)) + d\pgsx(v_i)\pfsx+ d\fsx \partial (\pgsx(v_i)).$$
Recall here that $f_{\si,x}$ was introduced in Proposition \ref{prop PMx}
and defined by the relation $\si = f_{\si,x}\exp(g_{x}) e^d$,
where the local section $\exp(g_{x}) e^d$ has norm one at $x$.
It is enough to bound by $O((\log d)^2)$ the integral 
$$ \frac{1}{d^k}\int_{\hld} || \log (\sumn |h_i|^2) (\ddbar \log (\sumn |h_i|^2))^{k}|| d\mu(\si),$$
since $||v_1||^2$ is bounded from below and above by positive constants in the model $(r)$.
Recall that 
$$ \ddbar (\sumn |h_i|^2) = \frac{\sumn \partial h_i\wedge \overline{\partial h_i}}{\sumn  |h_i|^2}+
  \frac{\sumn  h_i \overline{\partial h_i}\wedge \sum_{j=0}^{n-1}  \overline {h_j} \partial h_j}   {(\sumn  |h_i|^2)^2}.$$
From this we deduce, for $k\in \unan$, the upper bound 
$$|| (\ddbar \log(\sumn |h_i|^2))^k|| \leq \frac{2^k}{(\sumn |h_i|^2)^k} \sum_{|I|= k, |J|=k} ||\partial h_I\wedge \partial h_J||,$$
where  
$I$ and $J$ are ordered sets of $k$ elements in $\{1, \cdots, n\}$ and 
$\partial h_I= \partial h_{i_1}\wedge \cdots\wedge \partial h_{i_k}$ if $I = (i_1, \cdots , i_k)$. 
Our integral gets then bounded from above by
$$ \frac{2^k}{d^k}\sum_{|I|= k, |J|=k} \int_{\hld} \frac{|\log (\sumn |h_i|^2)|}{(\sumn |h_i|^2)^k} ||\partial h_I\wedge \partial h_J||
d\mu(\si).$$
Denote by $H^\perp_{x}$ the orthogonal complement of $H_{x}$ in $H^0(X;L^d)$, see \S \ref{II 3.2 jets}.
Likewise, with a slight abuse of notation, denote  by $H_{x}/H_{2x}$ (resp. $H_{2x}/H_{3x}$) the orthogonal complement 
of $H_{2x}$ (resp. $H_{3x}$) in $H_{x}$ (resp. $H_{2x})$. The space $H^0(X;L^d) $
then writes as a product 
$$H^0(X;L^d)=H^\perp_{x  }\times (H_{x}/H_{2x})\times (H_{2x}/H_{3x})\times H_{3x}$$
while its Gaussian measure $\mu$ is a product measure. The terms in our integral  only involve 
 jets at the second order of  sections and hence are constant on $H_{3x}$. Using Fubini's theorem,
it becomes thus enough to bound the integral over the space   $H^\perp_{x  }\times (H_{x}/H_{2x})\times (H_{2x}/H_{3x})$,
whose dimension no more depends on $d$. 

Moreover, the subspace $V^\perp\subset T_{\xz}^*X\otimes L^d_{\xz}$ of forms that vanish on the $n-1$
vectors $v_1(x), \cdots, v_{n-1}(x) $ given by the relative trivialization  is one-dimensional. It induces an 
orthogonal decomposition $T_{\xz}^*X\otimes L^d_{\xz} = V\oplus V^\perp$, where
$V$ is of dimension $n-1$. The inverse image of $V^\perp$ in $H_{x}/H_{2x}$ by the evaluation
map $eval_{2x}$ is the line $D$ of $H_{x}/H_{2x}$ containing the sections $\si$ of
$H_{x}$ whose derivatives at $x$ vanish against $v_1(x), \cdots, v_{n-1}(x) $. 
We denote by $\widetilde H_{x} $ the orthogonal complement of $D$ in  $H_{x}/H_{2x}$
and by $\widetilde H_{2x}$ the direct sum $D\oplus (H_{2x}/H_{3x})$.

We then write  $\si = (\si_0, \si_1, \si_2)\in H^\perp_{\xz}\times \widetilde H_{\xz}\times \widetilde H_{2\xz}$
 and   $|h_0(x)| = d\sqrt{h_d(\si_0,\si_0)}_{|x}=c_d ||\si_0||$, where $||\si_0||= \sqrt {\langle \si_0,\si_0\rangle}$
and 
by Proposition \ref{prop norm}, $c_d d^{-\frac{n+2}{2}}$
remains bounded  between positive constants as $d$ grows to infinity. For $i\in \unan$, $h_i$ 
linearly depends  on $\si_0$ and $\si_1$; we write $h_i(\si_0,\si_1)$ this linear expression. 
The derivatives $\partial h_I$ and $\overline{\partial h_J}$ depend on $\si_2$. The expression
$\partial h_I\wedge \overline{\partial h_J}$ expands as a sum of $3^k $ terms,  some 
of which vanish since
the forms $\pfsx$ and $\overline{\pfsx}$ can only appear once in the expression. Denote by $h_{I\bar{ J}}$
one of these $3^k$ terms. It is a monomial of degree $2k$ in $\si_0, \si_1, \si_2$
and we denote by $l_0$ the degree of $\si_0$, by 
$l_1$ the degree of $\si_1$ and by $l_2 =2k-l_0-l_1$ the degree of $\si_2$
in this monomial.
Now, it is enough to bound from above by a $O((\log d)^2)$ 
the following integral, where $I,J\subset \{1, \cdots , n\}^k$
are given:
\beq
&& \frac{1}{d^k}\int_{H^\perp_{x}\times\widetilde H_{x} \times \widetilde H_{2x}} 
\frac{|\log \sumn |h_i(x)|^2|}{ (\sumn |h_i(x)|^2)^k}|| h_{I\bar{ J}} (\si_0, \si_{1}, \si_{1})|| d\mu(\si)
\eeq
\beq
\leq 
\frac{1}{c_d^{2k}d^k}\int_{H^\perp_{x}\times\widetilde H_{x} \times \widetilde H_{2x}} 
&& \frac{|\log (c_d^2||\si_0||^2) + \log (1+ \sumunn \frac{1}{c_d^2}|h_i(\frac{\si_0}{||\si_0||},\frac{\si_1}{||\si_0||})|^2)|}
{(1+ \sumunn \frac{1}{c_d^2} |h_i(\frac{\si_0}{||\si_0||}, \frac{\si_1}{||\si_0||})|^2)^k||\si_0||^ {2k}}... \\
&...& ||h_{I\bar{ J}} (\frac{\si_0}{||\si_0||},\si_1,\si_2) || ||\si_0||^{l_0} 
e^{-||\si_0||^2 - ||\si_1||^2 } d\si_0 d\si_1 d\mu(\si_2). 
\eeq
We  replace $\frac{\si_0}{||\si_0||}$ by 1 in this integral, without loss of generality, since it remains bounded.
Define $\alpha_1  = \si_1/||\si_0||$, so that $d\alpha_1= \frac{1}{||\si_0||^{2(n-1)}} d\si_1$. The integral rewrites
\beq
\frac{1}{c_d^{2k}d^k}\int_{H^\perp_{x}\times\widetilde H_{x} \times \widetilde H_{2x}} 
&&\frac{|\log (c_d^2||\si_0||^2) + \log (1+ \sumunn \frac{1}{c_d^2}|h_i(1,\alpha_1)|^2)|}
{(1+ \sumunn \frac{1}{c_d^2} |h_i(1,\alpha_1)|^2)^k}... \\
&... &||h_{I\bar{ J}} (1,\alpha_1,\si_2) || ||\si_0||^{2(n-1)-l_2} 
e^{-||\si_0||^2(1+||\alpha_1||^2) } d\si_0 d\alpha_1d\mu(\si_2).
\eeq
Now set $\beta_0 = \si_0\sqrt{1+ ||\alpha_1||^2}$ and $\beta_1 = \frac{\alpha_1}{\sqrt d}$, so that 
$d\beta_0 d\beta_1 = \frac{1+||\alpha_1||^2}{d^{n-1} } d\si_0 d\alpha_1$. The integral becomes
\beq 
\frac{1}{c_d^{2k}d^{1+\frac{l_0}{2}}}\int_{H^\perp_{x}\times\widetilde H_{x} \times \widetilde H_{2x}} 
&\frac{\left|\log (c_d^2||\beta_0||^2) - \log (1+d||\beta_1||^2) + 
\log (1+ \sumunn \frac{1}{c_d^2} |h_i(1,\sqrt d \beta_1)|^2)\right|}
{(1+ \sumunn \frac{1}{c_d^2} |h_i(1,\sqrt d\beta_1)|^2)^k(\frac{1}{d} + ||\beta_1||^2) ^ {n-\frac{l_2}{2}}}... \\ 
...&  \ ||h_{I\bar{ J}} (1,\beta_1,\si_2) || ||\beta_0||^{2(n-1)-l_2} 
e^{-||\beta_0||^2 } d\beta_0 d\beta_1d\mu(\si_2).
\eeq
The only terms depending on $\beta_0$ in this integral are $\log ||\beta_0||^2$, $||\beta_0||^{2(n-1)-l_2}$
and $e^{-||\beta_0||^2}$. They can be extracted from it and integrated over $H_{x}^\perp$ thanks to  Fubini's theorem. 
The latter integral over $H_{x}^\perp$ turns out to be bounded independently of $d$. 
As a consequence, it becomes enough to bound from above the following 
\beq
\frac{1}{c_d^{2k}d^{1+\frac{l_0}{2}}}\int_{\widetilde H_{x}\times  \widetilde H_{2x}}  
&\frac{\left|\log (c_d^2/d) - \log (\frac{1}{d} +||\beta_1||^2) + 
\log (1+ \sumunn \frac{1}{c_d^2} |h_i(1,\sqrt d\beta_1)|^2)\right|}
{(1+ \sumunn \frac{1}{c_d^2} |h_i(1,\sqrt d\beta_1)|^2)^k(\frac{1}{d} + ||\beta_1||^2) ^ {n-\frac{l_2}{2}}}...\\ 
...& ||h_{I\bar{ J}} (1,\beta_1,\si_2) || d\beta_1d\mu(\si_2).
\eeq
Note that by Proposition \ref{prop norm} and by definition of the functions $h_i$, $1\leq i \leq n-1$,
the expressions $\frac{1}{c_d} h_i(1,\sqrt d \beta_1)$ are affine in $\beta_1$ with  coefficients bounded
independently of $d$. Indeed, only the first term 
$\partial f_{\si, x} (v_i)$ of $h_i = \pfsx(v_i) + d\fsx \pgsx(v_i)$
depends on $\beta_1$ and by Proposition \ref{prop norm},
$ ||\pfsx||^2_{|x} = \sqrt{h^d(\nabla \si_1,\nabla \si_1)}_{|x}$ grows
as $||\si_1||$ times $d^{\frac{n+1}{2}}$ while $c_d$ grows as  $d^{\frac{n+2}{2}}$.  
Let us denote by $g_i(\beta_1)$ these expressions.

Likewise, the monomial $h_{I\bar{ J}} $ is the product of three monomials of degrees $l_0$, $l_1$ and $l_2$
in $\si_0$, $\si_1$ and $\si_2$ respectively. By Proposition \ref{prop bounded}, the coefficients of these
monomials are $O(d^{\frac{n+2}{2}l_0})$,  $O(d^{\frac{n+1}{2}l_1+l'_1})$ and  $O(d^{\frac{n+2}{2}l_2})$ 
respectively, where $l'_1$ equals $0$ if neither $\pfsx$ nor $\overline{\pfsx}$  appears
 in the monomial
$h_{I\bar{ J}}$,  equals 
$1$ if one of these two forms appears and  2 if they  both appear. As a consequence, 
$||h_{I\bar{ J}} (1,\beta_1, \si_2)||$ is bounded from above,
up to a constant, by $d^{(n+2)k + l'_1 - l_1/2}||\beta_1||^{l_1} ||\si_2||^{l_2}$.
Now, only the term in $||\si_2||^{l_2} $ depends on $\si_2$ in our integral.
 Again using Fubini's theorem, we may 
first integrate  over $\widetilde H_{2x}$
 equipped with the Gaussian measure $d\mu(\si_2)$ to get an integral  bounded independently of $d$. 
The upshot is that we just need to bound the integral
$$I=\frac{1}{d^{1+\frac{l_0+l_1}{2}-l'_1}} \int_{\widetilde H_{x}}  
\frac{|\log (c_d^2/d) - \log (\frac{1}{d} +||\beta_1||^2) + 
\log (1+ \sumunn |g_i( \beta_1)|^2)|}
{(1+ \sumunn  |g_i( \beta_1)|^2)^k(\frac{1}{d} + ||\beta_1||^2) ^ {n-\frac{l_2}{2}}} 
||\beta_1||^{l_1} d\beta_1.$$
There is a compact subset $Q$ of $\widetilde H_{x}$ independent of $d$ and a constant
$C>0$ independent of $x$ and $d$ such that 
$$ \forall \beta_1 \in Q, 1\leq 1+ \sumunn  |g_i( \beta_1)|^2 \leq C \text { and }$$
$$  \forall \beta_1 \in \widetilde H_{x} \setminus  Q,  1+ \sumunn  |g_i( \beta_1)|^2 \geq \frac{1}{C}||\beta_1||^2,$$
since by Definition \ref{def adapted} the vector fields $v_i$ remain uniformly linearly independent
on the whole $U$. 

Bounding from above the term $||\beta_1||^{l_1}$ by $(\frac{1}{d} + ||\beta_1||^2)^{l_1/2}$, 
we finally just have to estimate  from above the integrals
$$I_1=  \frac{\log d}{d^{1+\frac{l_0+l_1}{2}-l'_1}} 
\int_Q \frac{1}{(\frac{1}{d} + ||\beta_1||^2)^{n-\frac{l_1+l_2 }{2}}}d\beta_1$$
and 
$$ I_2= \frac{1}{d^{1+\frac{l_0+l_1}{2}-l'_1}}\int_{\widetilde H_{x}\setminus Q} 
\frac{ \log d + \log (||\beta_1||^2) d\beta_1}{||\beta_1||^{2k}(\frac{1}{d} + ||\beta_1||^2)^{n-\frac{l_1+l_2 }{2}}},$$
since $\log (1/d + ||\beta_1||^2)$ over $Q$ and $\log (c_d/d) $ are $O(\log d). $
Note that  $l'_1\leq \max(2,l_1)$, so that the exponent $1+ \frac{l_0+l_1}{2}-l'_1$ is never negative and vanishes
if and only if $l_0 = 0$, $l_1=l'_1=2$ and thus $l_2 = 2k-2$. 
There exists $R>0$ such that 
$$  I_1 \leq \frac{Vol(S^{2n-3})\log d}{d^{1+\frac{l_0+l_1}{2}-l'_1}} \int_{\frac{1}{d}}^{R} \frac{du }{u^{2-\frac{l_1+l_2 }{2}}}du
= O((\log d)^2),
$$
where $u=\frac{1}{d} + ||\beta_1||^2$.
Likewise, there exists $T>0$ such that 
\beq
 \int_{\widetilde H_{x}\setminus Q} 
\frac{\log d +\log ||\beta_1||^2 d\beta_1}{||\beta_1||^{2k}(\frac{1}{d} + ||\beta_1||^2)^{n-\frac{l_1+l_2 }{2}}}
&\leq & Vol(S^{2n-3}) \int_{T}^{\infty} \frac{ \log d+\log u }{u^{2+l_0/2}} du.
\eeq
This last integral is a $O(\log d)$, which implies Proposition \ref{prop bounded} when $U$ is of type $(r)$.

\subsubsection{Proof of Proposition \ref{prop bounded} along the base and critical loci of $p$}\label{para b}
The compact $K$ given by Proposition \ref{prop bounded} is covered by a finite number
of elements of the atlas $\mathcal U$ adapted to $(p,L)$. Moreover, we may assume that such elements are
all disjoint from the critical locus $\text{Crit}(p)$ when $n>2$. Proposition \ref{prop bounded} was proved
in \S \ref{para 4.1} for elements $U$ of type $(r)$ of $\mathcal U$. Let us assume now that $U$ is such an element of
type $(b)$ given by Definition \ref{def adapted} and that  $x\in K \cap U$. The case where $U$ is of 
type $(c)$ when $n=2$ just follows along the same lines. The main part of the proof
is similar to the one given in \S \ref{para 4.1} and we have to bound from above the integral
$$I = \frac{1}{d^{1+\frac{l_0+l_1}{2}-l'_1}} \int_{\widetilde H_{x}}  
\frac{|\log (c_d^2/d) - \log (\frac{1}{d} +||\beta_1||^2) + 
\log (1+ \sumunn |g_i( \beta_1)|^2)|}
{(1+ \sumunn  |g_i( \beta_1)|^2)^k(\frac{1}{d} + ||\beta_1||^2) ^ {n-\frac{l_2}{2}}} 
||\beta_1||^{l_1} d\beta_1.$$
However, the norm $||v_1||$ of the vector $v_1$ given by Definition \ref{def trivialization}  
converges now to $0$ when $x$ approaches the base locus of $p$. 
Denote  by $\widetilde H_{x}''$ the hyperplane of $\widetilde H_{x}$ consisting of the sections
whose 1-jet at $x$ vanish against $v_1$, that is the sections whose image under
 $eval_{2x}$ vanish against $v_1$, see \S \ref{II 3.2 jets}. 
Denote then by  $\widetilde H_{x}'$ the line  orthogonal to $\widetilde H_{x}''$ in  $\widetilde H_{x}$
and by $\beta_1 = (\beta_1', \beta_1'') $ the coordinates on $\widetilde H_{x}=\widetilde H_{x}'\times \widetilde H_{x}''$. 
This time there exists a compact subset $Q =Q'\times Q''$ of $\widetilde H_{x}$, independent of $d$ and of $x\in K\cap U$,
as well as a constant $C>0$ such that 
$$ \forall \beta_1 \in Q, 1\leq 1+ \sumunn  |g_i( \beta_1)|^2 \leq C \text { and }$$
$$  \forall \beta_1 \in \widetilde H_{x} \setminus  Q,  1+ \sumunn  |g_i( \beta_1)|^2 \geq \frac{1}{C} 
(1+ ||v_1(x)||^2 ||\beta_1'||^2 + ||\beta''_1||^2)
.$$
The integral $I$ over the compact $Q$ is bounded from above by a $O((\log d)^2)$, see \S \ref{para 4.1}.
Only the second integral differs. In order to estimate the latter, let us bound from above $||\beta_1||^{l_1}$ by 
$ (\frac{1}{d} + ||\beta_1||^2)^{l_1/2}$. We have  to bound  the integral
$$\frac{1}{d^{1+\frac{l_0+l_1}{2}-l'_1}} \int_{\widetilde H_{x} \setminus  Q}  
\frac{\log d + \log ||\beta_1||^2 d\beta_1}
{(1+ ||v_1||^2 ||\beta_1'||^2 + ||\beta''_1||^2)^k(\frac{1}{d}+  ||\beta_1||^2)^{n- \frac{l_1+l_2}{2}}}
.$$
When $n- \frac{l_1+l_2}{2} >1$, let us bound from above this integral  by 
$$ \frac{1}{d^{1+\frac{l_0+l_1}{2}-l'_1}} \int_{\widetilde H''_{x} \setminus  Q''} 
\frac{d\beta_1''}{||\beta_1''||^{2k}}
 \int_{\widetilde H'_{x} \setminus  Q'} 
\frac{\log d + \log ||\beta_1||^2  d\beta'_1}
{( ||\beta_1||^2)^{n- \frac{l_1+l_2}{2}}}
.$$
There exists $ R>0$ such that 
\beq
 \frac{1}{\pi}  \int_{\widetilde H'_{x} \setminus  Q'} 
\frac{\log d + \log ||\beta_1||^2 d\beta'_1}
{( ||\beta_1||^2)^{n- \frac{l_1+l_2}{2}}}
&\leq& \int_R^\infty 
\frac{\log d +\log ( ||\beta''_1||^2 + u)du}
{(  ||\beta''_1||^2+ u)^{n- \frac{l_1+l_2}{2}}}\\
&=& \left[ \frac{\log d+\log ( ||\beta''_1||^2 + u)}
{   ( \frac{l_1+l_2}{2}-n+1)(||\beta''_1||^2+ u)^{n- \frac{l_1+l_2}{2}-1}  }      \right]^\infty_R\\
&&+  \int_R^\infty \frac{du}
{(n- \frac{l_1+l_2}{2}-1)(||\beta''_1||^2+ u)^{n- \frac{l_1+l_2}{2} }   }  \\
&=& 
\frac{\log d + \log (||\beta''_1||^2 + R)}
{( n-\frac{l_1+l_2}{2}-1)(  ||\beta''_1||^2+ R)^{n- \frac{l_1+l_2}{2}-1}      }    \\
& & + \frac{1}{(n- \frac{l_1+l_2}{2}-1)^2(  ||\beta''_1||^2+ R)^{n- \frac{l_1+l_2}{2}-1}   }
\eeq
Hence, our integral gets bounded from above,  up to a constant, by the integral 
$$ 
\int_{\widetilde H''_{x} \setminus  Q''} \frac{\log d + \log (||\beta_1''||^2+R  ) }{ ||\beta_1''||^{2(n-1)+l_0} } d\beta_1'',$$
which is itself a $O(\log d)$ since $\widetilde H''_{x}$ is of dimension $n-2$. 
When $n- \frac{l_1+l_2}{2} =1$, which implies that $l_0=0$ and $k=n-1$, we observe that
\beq
 (1+ ||v_1||^2 ||\beta_1'||^2 + ||\beta''_1||^2)^k &=& (1+ ||v_1||^2 ||\beta_1'||^2 + ||\beta''_1||^2) (1+ ||v_1||^2 ||\beta_1'||^2 + ||\beta''_1||^2)^{k-1}\\
& \geq & ||v_1||^2 (1+||\beta_1||^2)(1+||\beta_1''||^2)^{k-1}\\
& \geq & ||v_1||^2 ||\beta_1||^2||\beta_1''||^{2(k-1)}
\eeq
as long as $||v_1||\leq 1$, which can be assumed. 
Our integral gets then bounded by 
$$ \frac{1}{d^{1+\frac{l_0+l_1}{2}-l'_1}||v_1||^2} 
\int_{\widetilde H''_{x} \setminus  Q''} 
\frac{d\beta_1''}{||\beta_1''||^{2(k-1)} }
 \int_{\widetilde H'_{x} \setminus  Q'} 
\frac{\log d + \log ||\beta_1||^2  d\beta'_1}
{  ||\beta_1||^4}
$$
which is similar to the previous one. The latter is then bounded from above
by a $O(\frac{\log d}{||v_1||^2})$, implying the result.

\section{Total Betti numbers of random real hypersurfaces}\label{section real}
\subsection{Statement of the results}\label{para 3.1}
\subsubsection{Expectation of the total Betti number of  real hypersurfaces}\label{para 1.1}
Let $(X,c_X)$ be a  smooth real projective manifold of positive dimension $n$, meaning that
$X$ is a smooth  $n$-dimensional complex projective  manifold equipped with an antiholomorphic
involution $c_X$. Let $\pi: (L, c_L) \to (X,c_X)$ be a real holomorphic ample  line bundle,
so that the antiholomorphic involutions satisfy $ \pi\circ c_L = c_X\circ \pi$. 
For every $d>0$, we denote by $L^d$  the $d$-th tensor power of $L$,
by $\Rr H^0(X;L^d)$ the space of  global real holomorphic sections  of $L^d$,
which are the sections $\si\in H^0(X;L^d)$ satisfying $\si\circ c_X = c_L\circ \si$,
and by $\Rr \Delta_d= \Delta_d \cap \Rr H^0(X;L^d)$ the real discriminant locus.

For every section $\si \in \Rr H^0(X;L^d)\setminus \Rr \Delta_d$, $C_\si = \si^{-1}(0)$
is a smooth real  hypersurface of $X$. By Smith-Thom's inequality, see Theorem \ref{theo Smith},
the total Betti number
 $b_*(\Rr C_\si;\zsz) = \sumn \dim H_i(\Rr C_\si, \Zz/2\Zz) $ of its real locus is bounded
 from above by the total Betti number $b_*(C_\si;\zsz) = \sum_{i=0}^{2n-2} \dim H_i(C_\si; \Zz/2\Zz) $ 
 of its complex locus, which  from Lemma \ref{lemma Betti asymptotic} 
is equivalent to $(\int_X c_1(L)^n) d^n$ as $d$
 grows infinity. What is the expectation of this real total
 Betti number? If we are not able to answer to this question,
we will estimate this number from above, see Theorems
\ref{theo principal} and \ref{theo principal produit}. Note
that in dimension one, such an upper bound can be deduced from our recent
work \cite{GW}. 
 
Let us first precise the measure of probability considered on 
 $\Rr H^0(X;L^d)$. We proceed as in \S \ref{section random divisors}.
 We equip $L$ with a real Hermitian metric $h$ of positive curvature $\omega \in \Omega ^{(1,1)}(X,\Rr)$,
real meaning that $c_L^*h = h$. As in \S \ref{para notations}, we denote by $dx = \frac{1}{\int_X \omega^n }\omega^n$
 the associated  volume form of $X$, by 
 $\langle \ \rangle$ the induced $L^2$-scalar product  on $\Rr H^0(X;L^d)$,
 and by $\mu_\Rr$ the associated Gaussian measure, defined by the relation 
 $$ \forall A \subset \Rr H^0(X;L^d), \mu_\Rr (A) = \frac{1}{(\sqrt \pi )^{N_d}  }\int_A e^{-||\si||^2}d\si.$$
For every $d>0$, we denote by 
$$E_\Rr(b_*(\Rr C_\si;\zsz )) = \int_{\Rr H^0(X;L^d)\setminus \Rr \Delta_d} b_*(\Rr C_\si;\zsz) d\mu_\Rr (\si)  $$
the expected total Betti number of  real hypersurfaces linearly equivalent to $L^d$.
\begin{theorem}\label{theo principal}
 Let $(X,c_X)$ be a smooth real projective manifold of dimension $n$ greater than one
 equipped with a Hermitian real line bundle $(L,c_L)$
of positive curvature.  Then, the expected total Betti number 
$E_\Rr(b_*(\Rr C_\si;\zsz))$ is a $o(d^{n})$ and even a $O(d(\log d)^2)$
if $n=2$.
\end{theorem}
Note that the exact value of the expectation
$E_\Rr(b_*(\Rr C_\si; \zsz))$ is only known  when $X=\Cc P^1$, see \cite{Kac}, \cite{SS} and \cite{EK}. 
While we were writing this article in june 2011, 
Peter Sarnak  informed us that together with Igor Wigman, he can bound this expectation by a $O(d)$
when $X= \Cc P^2$ and suspects it is equivalent to a constant times $d$ when $d$ grows 
to infinity.  Such a guess was already made couple of years ago by Christophe Raffalli, 
 based on computer experiments. 

It could be that
this expectation is in fact equivalent to $d^{\frac{n}{2}}$ times a constant as soon as the real locus of the 
manifold $(X,c_X)$ is non empty and 
in particular that the bound given
by Theorem \ref{theo principal} can be improved by a $O(d^{\frac{n}{2}})$. 
When $(X,c_X)$ is a product of smooth real projective curve for instance, we can improve the 
$o(d^n)$ given by Theorem \ref{theo principal} by a $O(d^{\frac{n}{2}}(\log d)^n)$,
being much closer to a $O(d^{\frac{n}{2}})$ bound, see Theorem \ref{theo principal produit} below.
This $d^{\frac{n}{2}}$
can be understood as the volume of a $\frac{1 }{\sqrt d}$ neighborhood of the real locus $\Rr X$ in $X$ for the volume
form induced by the curvature of $L^d$, where $\frac{1 }{\sqrt d}$ is a fundamental scale in \kah geometry and 
H\"ormander's theory of peak sections. A peak section 
centered at $x$ can be symmetrized to provide a real  section having two peaks near $x$ and $c_X(x)$, 
see \S\ref{para real H}.
This phenomenon plays an important r\^ole in the proof of Theorem \ref{theo principal}  and 
seems to be intimately related to the value of the expectation $E_\Rr(b_*(\Rr C_\si;\zsz))$.  
Note finally 
that Theorem \ref{theo principal} contrasts with the computations made by Fedor Nazarov and
Mikhail Sodin in \cite{Nazarov-Sodin} for spherical harmonics in dimension two, as well as with the one achieved
by Maria Nastasescu for a real Fubini-Study measure, as P. Sarnak informed us. In both cases, the expectation is quadratic.
\begin{theorem}\label{theo principal produit}
Let $(X,c_X)$ be the product of $n>1$ smooth real projective curves, equipped with 
a real Hermitian line bundle  $(L,c_L)$ of positive curvature. Then, the expected
 total Betti number  $E_\Rr (b_*(\Rr C_\si ; \zsz)) $ is a $O(d^{\frac{n}{2}} (\log d)^n)$. 
\end{theorem}

\subsubsection{Random real divisors and distribution of critical points}
Our proof of Theorem \ref{theo principal} is based on a real analogue of Theorem \ref{theo equi}
that we formulate here, see Theorem \ref{theo equi real}. We use the notations introduced in \S \ref{para 1.1} and equip $X$ 
with a real Lefschetz pencil $p: X\dashrightarrow \Cc P^1$, see \S \ref{section Betti}. 
For every $d>0$, denote by $\Rr \widetilde \Delta_d = \widetilde \Delta_d \cap \Rr H^0(X;L^d)$ the 
union of $\Rr \Delta_d$ with the set of sections $\si \in \Rr H^0(X;L^d)$ 
such that  either $C_\si $ contains a critical point of $p$,
or $p_{|C_\si}  $ is not Lefschetz, see \S \ref{para notations}. Denote, as in \S \ref{para notations},
by $\mathcal R_\si$ the  critical locus  of the restriction $p_{|C_\si}$, 
where $\si \in \rhld$. Then, for every continuous function $\chi: X \to \Rr$,
 denote by $$E_\Rr (\langle \nu_\si, \chi \rangle) = \frac{1}{\#\mathcal R_\si}\int_{\rhld} (\sum_{x\in \mathcal R_\si}\chi(x))
d\mu_\Rr(\si)$$
the expectation of the probability measure $\nu_\si$ carried by the critical points of $p_{|C_\si}$,
see Definition \ref{def Dirac}, computed with respect to the real Gaussian measure $\mu_\Rr$ and evaluated against $\chi$.
\begin{theorem}\label{theo equi real}
 Let $(X,c_X)$ be a smooth real projective  manifold of positive dimension $n$ equipped with a real Lefschetz pencil 
$p: X\dashrightarrow \Cc P^1$  of critical locus $\text{Crit}(p)$. Let $(L,c_L) \to (X,c_X)$ 
be a real ample line bundle equipped with a real Hermitian metric of positive curvature. Let $(\chi_d)_{d\in \Nn^*}$
be a sequence of elements of $C^2(X,\Rr)$ which converges to $\chi$ in $L^1(X,\Rr)$ as 
$d$ grows to infinity. Assume that
$$\sup_{x\in \text{Supp}(\chi_d)} d(x,c_X(x))>  2\frac{\log d }{\sqrt d}.$$
When $n>2$, assume moreover that
the distance between $\text{Crit}(p)$ and the supports of $\chi_d$, $d>0$, are uniformly bounded from below
by some positive constant.
Then, the real expectation $E_\Rr (\langle \nu_\si, \chi \rangle) $ converges to the integral $\int_X \chi dx$ 
as $d$ grows to infinity. More precisely,
$$  E_\Rr (\langle \nu_\si, \chi \rangle) = \int_X \chi_d dx + 
O\left(\frac{(\log d)^2}{d} ||\ddbar \chi_d||_{L^1} \right)+ O\left(\frac{1}{d} ||\chi_d||_{L^1}\right).$$
\end{theorem}
In this Theorem \ref{theo equi real}, $\text{Supp}(\chi_d)$ denotes the support of $ \chi_d$
and $d(x,c_X(x))$  the distance between the points $x$ and $c_X(x)$ for the \kah metric 
induced by the curvature $\omega$ of $L$.

\subsection{Real peak sections and evaluation of two-jets of sections}\label{para real H}
Let $(X,c_X)$ be a smooth real projective manifold of positive dimension $n$ and
$(L,c_L) \to (X,c_X)$ be a real holomorphic ample line bundle 
equipped with a real Hermitian metric $h$ of positive curvature.
\begin{definition}\label{def peak real}
A real peak section of the ample real holomorphic line bundle $(L,c_L)$
over the real projective manifold $(X,c_X)$  is a section which  writes
 $\frac{\si + c^*\si}{||\si + c^*\si ||}$, where $\si$ is a  
peak section of \hor given by Definition \ref{def peak} and  $c^*\si = c_{L^d} \circ \si \circ c_X$.
\end{definition}
Recall that from  Lemma \ref{lemma peak}, the $L^2$-norm of a peak section
concentrates in a $\frac{\log d }{\sqrt d}$ neighborhood of a point $x\in X$. 
When $x$ is real, the real peak section $\frac{\si + c^*\si}{||\si + c^*\si||}$
looks like a section of \hor given by Definition \ref{def peak}.
When the distance between $x$ and $c_X(x)$ is bigger than $\frac{\log d }{\sqrt d}$, more or less half of 
the $L^2$-norm of this real section concentrates in a neighborhood of $x$ 
and another half in a neighborhood of $c_X(x)$. 
Such a real section has thus 
two peaks near $x$ and $c_X(x)$. When $d(x,c_X(x))< 2\frac{\log d }{\sqrt d}$,
these two peaks interfere, interpolating the  extreme cases  just
discussed. We are now interested in the case  $d(x,c_X(x))> 2\frac{\log d }{\sqrt d}$,
where we can establish a real analogue of Proposition \ref{prop norm}.
\begin{lemma}\label{lemma real scal}
Let $(L,c_L)$ be a real holomorphic Hermitian line bundle of positive curvature over a smooth real  projective
manifold $(X,c_X)$ of positive dimension $n$. Let $(x_d)_{d\in \Nn^*}$ be a sequence of points
such that $ d(x_d, c_X(x_d))>2\frac{\log d }{\sqrt d}$  and 
let $(\si_d)_{d\in \Nn^*}$ be an associated
sequence of sections given by Lemma \ref{lemma peak} with $p'=2$. 
Then, the Hermitian product $\langle \si_d,c^*\si_d \rangle$ is a $O(\frac{1}{d^2})$, so that 
the norm $||\si_d+  c^*\si_d||$ equals $\sqrt 2 ||\si_d||(1 + O(1/d^2))$. 
\end{lemma}
\bpr 
By definition,
\beq 
\langle \si_d,c^*\si_d \rangle &=& \int_X h^d(\si_d,c^*\si_d) dx\\
&=&  \int_{B(\xz, \frac{\log d }{\sqrt d})} h^d(\si_d,c^*\si_d) dx + \int_{X\setminus B(\xz, \frac{\log d }{\sqrt d})} 
h^d(\si_d,c^*\si_d) dx\\
& \leq & \left(\int_{B(\xz, \frac{\log d }{\sqrt d})} h^d(\si_d, \si_d) dx\right) ^{1/2}\left(\int_{B(\xz, \frac{\log d }{\sqrt d})} 
h^d(c^*\si_d,c^*\si_d)dx\right) ^{1/2}\\
&&+ \left(\int_{X\setminus B(\xz, \frac{\log d }{\sqrt d})} h^d(\si_d, \si_d) dx\right) ^{1/2}\left( \int_{X\setminus B(\xz, \frac{\log d }{\sqrt d})} 
h^d(c^*\si_d,c^*\si_d)dx\right) ^{1/2} 
\eeq
by Cauchy-Schwarz's inequality. We deduce that
\beq
\langle \si_d,c^*\si_d \rangle &\leq & ||\si_d|| \left[\left( \int_{X\setminus B(\xz,\frac{\log d}{\sqrt d})} h^d(\si_d, \si_d) dx\right) ^{1/2} +
\left(\int_{B(\xz, \frac{\log d}{\sqrt d})} 
h^d(c^*\si_d,c^*\si_d)dx\right) ^{1/2}\right]
\eeq
By assumption, the balls $B(x,\frac{\log d }{\sqrt d})$ and 
$B(c_X(x),\frac{\log d }{\sqrt d})$ are disjoint, so that by 
 Lemma \ref{lemma peak}, these two last terms are $O(1/d^2)||\si_d||$. 
Hence, 
\beq 
\langle \si_d+ c^*\si_d, \si_d+ c^*\si_d\rangle &=& 2 ||\si_d||^2 + 2\Re \langle \si_d, 
c^*\si_d \rangle\\
& = & 2 ||\si_d||^2 + O(1/d^2) ||\si_d||^2,
\eeq
so that $||\si_d+ c^*\si_d|| = \sqrt 2 ||\si_d||(1 + O(1/d^2))$.
\epr
Set $\rhx = H_{\xz}\cap \rhl$, $\rhdx = H_{2\xz}\cap \rhl$ and $\rhtx = H_{3\xz}\cap \rhl$, 
where  $H_{\xz}$, $ H_{2\xz}$ and $H_{3\xz}$ have been introduced in \S \ref{II 3.2 jets}.
Likewise, with a slight abuse of notation, denote by $eval_{\xz}$, $eval_{2\xz}$ and $eval_{3\xz}$ the restrictions
of the evaluation maps to the spaces $\rhl/\Rr H_{\xz}$, $\Rr H_{\xz}/\Rr H_{2\xz}$ and $ \Rr H_{2\xz}/\Rr H_{3\xz}$
respectively, so that:
\beq
eval_{x}: \si \in \Rr H^0(X;L^d)/\Rr H_{x} & \mapsto &\si(x) \in L^d_{x},\\
eval_{2x}: \si \in \Rr H_{x}/\Rr H_{2x} &\mapsto &\nabla \si_{|x} \in T^*_{x}X\otimes L^d_{x},\\
\text{and }  eval_{3x}: \si \in \Rr H_{2x}/\Rr H_{3_{x}} & \mapsto &
\nabla^2\si_{|x} \in Sym^2 (T^*_{x}X)\otimes L^d_{x}.
\eeq
The following Proposition \ref{prop real norm} is a real analogue of Proposition \ref{prop norm}.
\begin{proposition}\label {prop real norm}
Let $(L,c_L)$ be a real holomorphic Hermitian line bundle of positive curvature over a smooth complex projective
manifold $X$ of positive dimension $n$. Let $(x_d)_{d\in \Nn^*}$ be a sequence of points in $X$
such that $ d(x_d, c_X(x_d))>2\log d /\sqrt d$.  
 Then, the maps
$d^{-\frac{n}{2}} eval_{x_d}$, $d^{-\frac{n+1}{2}} eval_{2x_d}$ and $ d^{-\frac{n+2}{2}} eval_{3x_d}$ as
well as their inverse have bounded norms and determinants, as
long as $d$ is large enough.
\end{proposition}
Note that the evaluation maps $eval_{x_d}$, $eval_{2x_d}$ and $eval_{3x_d}$
of Proposition \ref{prop real norm} are only $\Rr$-linear.
\bpr[ of Proposition \ref{prop real norm}]The proof is analogous to the one of Lemma \ref{lemma real scal}. Let 
$\si_0 = \frac{\si+ c^*\si}{||\si+ c^*\si||}$ be
a real peak section given by Definition \ref{def peak real}, where $\si$ 
is a  peak section given by $x_d$ and Definition \ref{def peak}. From Lemma \ref{lemma real scal},
we know that $||\si+ c^*\si||$ equals $\sqrt 2  + O(1/d^2)$. 
Moreover $2\frac{\log d }{\sqrt d}< d(x_d,c_X(x_d))$ and by Lemma \ref{lemma peak},
the $L^2$-norm of $c^*\si$  in a neighborhood of $x_d$ is a $O(1/d^2)$.
From the mean inequality, we deduce the bound (see \cite{Hormander}, Theorem 4.2.13 for instance)
$ h^d(c^*\si,  c^*\si)_{|x_d} = O(1/d^2),$
so that by Lemma 2.1 of \cite{Tian}, $d^{-n}h^d(\si_0,\si_0)_{|x_d}$ converges 
to a non negative constant as $d$ grows to infinity. Denote by $\si_0^{\rhx}$ 
the orthogonal projection of $\si_0$ onto $\rhx$. We proceed as in the proof of Proposition 
\ref{prop norm} 
to get that $\langle \si , \frac{\si_0^{\rhx}}{||\si_0^{\rhx}||}\rangle = O(\frac{1}{d})$
and  likewise $ \langle c^*\si , \frac{\si_0^{\rhx}}{||\si_0^{\rhx}||}\rangle = O(\frac{1}{d})$,
since $c$ is an isometry for the $L^2$-Hermitian product. 
Hence, $\langle \si_0 , \frac{\si_0^{\rhx}}{||\si_0^{\rhx}||}\rangle=O(\frac{1}{d})$. Writing
$ \si_0^\perp = \frac{\si_0 - \si_0^{\rhx}}{||\si_0 - \si_0^{\rhx}||},$ 
we deduce as in the proof of Proposition \ref{prop norm} 
that $d^{-n} h^d(\si_0^\perp, \si_0^\perp)_{|x_d}$ converges to a positive constant 
as $d$ grows to infinity. Replacing $\si$ by $i\si$, we
define 
$ \widetilde \si_0 =  \frac{i\si_0 + c^*(i\si_0) }{||i\si_0 + c^*(i\si_0)||}$
and check in the same way that $d^{-n} h^d(\widetilde \si_0^\perp, \widetilde \si_0^\perp)_{|x_d}$ converges to a positive constant 
as $d$ grows to infinity. But the quotient $\frac{\widetilde \si_0(x_d)}{\si_0(x_d)}$
and thus $\frac{\widetilde \si^\perp_0(x_d)}{\si_0^\perp(x_d)}$ converge to $i$
as $d$ grows to infinity. Likewise,
\beq
 \langle \si_0, \widetilde \si_0 \rangle&=& 
\frac{\langle \si+ c^*\si,  i\si-ic^*\si\rangle }
{||\si+ c^*\si||||i\si-i c^*\si|| }
= \frac{2\Re (i \langle \si, c^*\si\rangle) }
{||\si+ c^*\si||||i\si-i c^*\si|| } = O(\frac{1}{d}),
\eeq
so that $\langle \si_0^\perp, \widetilde \si_0^\perp\rangle = O(\frac{1}{d})$ and $d^{-\frac{n}{2}}eval_{x}$
as well as its inverse have bounded norms and determinants when $d$ grows to infinity. 
The remaining cases are obtained by similar modifications of the proof of 
Corollary \ref{cor expectation}.
\epr
\subsection{Proof of the main results}
\subsubsection{Proof of Theorem \ref{theo equi real}}
The proof goes along the same lines as  the one of Theorem \ref{theo equi}. We begin
with the following analogue of Corollary  \ref{cor expectation}.
\begin{corollary}\label{cor real expectation}
 Under the hypotheses of Proposition \ref{prop PMx}, we assume moreover
that the manifold $X$ and  Hermitian bundle $L$ are real. Then, 
\beq
 E_\Rr (\langle \nu_\si, \chi \rangle) &= &\frac{1}{\# \br_\si}d^n
\int_X  \chi  \om^n  \\ 
+& &\frac{1}{\# \br_\si} \sum_{k=0}^{n-1} (\frac{i}{2\pi})^{n-k}d^k 
\int_X  \ddbar \chi \wedge \omega^k\wedge 
\int_{\Rr H^0(X;L^d)\setminus \Rr \widetilde \Delta_d}\lambda_x(\ddbar \lambda_x)^{n-1-k} d\mu(\si),
\eeq
where $$\lambda_x =  
 \log \left(d^2 |f_{\si,x}|^2 + \sum^{n-1}_{i=1} |\partial f_{\si,x}(v_i)+f_{\si,x} \partial g_x(v_i)|^2\right).$$
\end{corollary}
\bpr
The result follows after integration over $\Rr H^0(X;L^d)\setminus \Rr \widetilde \Delta_d$ of the relation given by Proposition \ref{prop PMx}. 
\epr
\bpr[ of Theorem \ref{theo equi real}] 
We apply  Corollary \ref{cor real expectation} to the functions $\chi_d$, $d>0$ and use the notations of this corollary. 
The first term in the right hand side of the formula given  by this corollary
is $ \int_X \chi_d dx + O(\frac{1}{d}||\chi_d||_{L^1})$ as follows from  Proposition \ref{prop critical}.
It is thus enough to prove that for $k\in \{0, \cdots , n-1\}$ the integral 
$\int_{\Rr H^0(X;L^d)\setminus  \Rr \widetilde \Delta_d}\lambda_x(\ddbar \lambda_x)^{k} d\mu(\si)$
is uniformly bounded on the support of $\ddbar \chi_d$ by a $O(d^k(\log d)^2/||v_1||^2)$ since
$1/||v_1||^2$ is integrable over $X$. The space  $\Rr H^0(X;L^d)\setminus \Rr \widetilde \Delta_d$ 
is equipped with its Gaussian measure, but
the integrand  only depends on  the two-jets of sections $\si \in \rhl$ at
the point $x\in X\setminus (\Rr X \cup \text{Crit}(p))$. Since the Gaussian measure is a product measure, 
writing $\rhl$ as the product of the space $\Rr H_{3x}$ introduced in \S \ref{para real H}
with its orthogonal complement $\Rr H_{3x}^\perp$, we deduce that it is enough to prove this uniform bound
for the integral $ \int_{\Rr H_{3x}^\perp\setminus  \Rr \widetilde \Delta_d}\lambda_x(\ddbar \lambda_x)^{k} d\mu(\si)$.
By Proposition \ref{prop real norm}, the evaluation maps $eval_x$, $eval_{2x}$
and $eval_{3x}$ of jets up to order two at $x$  provide an isomorphism between 
$\Rr H_{3x}^\perp$ and 
$L_x^d\oplus (T^*_xX\otimes L^d_x) \oplus (Sym^2(T^*_xX)\otimes L^d_x)$. 
In particular, this implies that the space $\Rr H_{3x}^\perp$ is also a complement of the 
space $H_{3x}$ introduced in \S \ref{II 3.2 jets}. By Proposition \ref{prop norm},
these evaluation maps factor through an isomorphism $I^{-1}_d$ between $\Rr H_{3x}^\perp$ 
and the orthogonal complement $H^\perp_{3x}$ of $H_{3x}$ in $H^0(X;L^d)$. 
By  Proposition \ref{prop norm} and Proposition \ref{prop real norm}, the Jacobian 
of $I_d$ is bounded independently of $d$, while the Gaussian measure 
$\mu_\Rr$ is bounded from above on $\Rr H_{3x}^\perp$ by the measure $(I_d)_* \mu$ 
up to a positive dilation of the norm, independent of $d$. 
 After a change
of variables given by this isomorphism $I_d$ and after the one given by the dilation, 
it becomes enough to prove that the integral
$ \int_{ H^0(X;L^d)\setminus   \widetilde \Delta_d}\lambda_x(\ddbar \lambda_x)^{k} d\mu(\si)$
is uniformly bounded on the support of $ \chi_d$ by a $O(d^k(\log d)^2/||v_1||^2)$. The latter follows
from Proposition \ref{prop bounded}.
\epr
\subsubsection{Proof of Theorem \ref{theo principal}}\label{para 3.3.2}
Equip $(X,c_X)$ with a real Lefschetz pencil $p : X  \dashrightarrow \cpun$ 
and denote by $F$ a regular fiber of $p$. For every section $\si \in \rhl\setminus \Rr \widetilde \Delta_d$,
the restriction of $p$ to $\Rr C_\si$ satisfies the hypotheses of Lemma \ref{lemma Betti},
so that 
$$b_*(\Rr C_\si ; \zsz)\leq  b_*(\Rr F\cap \Rr C_\si ; \zsz) + \# \text{Crit}(p_{|\Rr C_\si}).$$ 
By  Smith-Thom's inequality, see Theorem \ref{theo Smith}, 
$b_*(\Rr F\cap \Rr C_\si ; \zsz)\leq b_*(F\cap C_\si ; \zsz)$
while from Lemma \ref{lemma Betti asymptotic} applied to $L_{|F}$, $b_*(F\cap C_\si ; \zsz)$
is a $O(d^{n-1})$. As a consequence, we have to prove that the expectation of $\#\text{Crit}(p_{|\Rr C_\si})$
is a $o(d^n)$ and even a $O(d(log d)^2)$ when $n=2$. 

Let us identify a neighborhood $V$ of $\Rr X$ in $X$ with the cotangent bundle of $\Rr X$.
We can assume that $V\setminus \Rr X$ does not contain any critical point of $p$.
Let $\chi : X \to [0,1]$ a be function of class $C^2$ satisfying $\chi =1$ outside of a
compact subset of $V$ and $\chi = 0$ in a neighborhood of $\Rr X$.
For every $d>0$, let $\chi_d : X \to [0,1]$ be the function which equals one 
outside of $V$ and whose restriction to $V$ writes in local coordinates 
$(q,p)\in V \simeq T^*\Rr X \mapsto \chi(q, \frac{\sqrt d}{\log d } p) \in [0,1]$,
where $q$ is the coordinate along $\Rr X$ and $p$ the coordinate along the
fibers of $T^*\Rr X$. This sequence $(\chi_d)_{d>0} $ converges to the constant function 1 
in $L^1(X,\Rr)$ as $d$ grows to infinity, while the norm $||\ddbar \chi_d||_{L^1(X,\Rr)}$
is a $O((\frac{\log d }{\sqrt d})^{n-2})$. Moreover, for every $x\in \text{Supp}(\chi_d)$,
$ d(x,c_X(x))>2 \frac{\log d }{\sqrt d}$, so that when $n=2$, 
Theorem \ref{theo equi real} applies. From Proposition \ref{prop critical}, 
we thus deduce that 
$$ \int_{\rhl \setminus \Rr \widetilde \Delta_d} 
(\sum_{x \in \mathcal R_\si} \chi_d (x)) d\mu_\Rr (\si) =\# \mathcal R_\si + O(d(\log d)^2).$$
Moreover, for every $\si \in \rhl \setminus \Rr \widetilde \Delta_d$ and every $d>0$,
we have
$ \# \text{Crit}(p_{|\Rr C_\si}) \leq \# \mathcal R_\si - \sum_{x \in \mathcal R_\si}\chi_d (x),$
so that after integration,
$$\int_{\rhl \setminus \Rr \widetilde \Delta_d} 
\# \text{Crit}(p_{|\Rr C_\si}) d\mu_\Rr (\si)=O(d(\log d)^2),$$
hence the result when $n=2$.

When $n>2$, we apply Theorem \ref{theo equi real} to the function $\chi$
and deduce that
$$ \int_{\rhl \setminus \Rr \widetilde \Delta_d} 
(\sum_{x \in \mathcal R_\si} \chi (x)) d\mu_\Rr (\si) = (\int_X \chi dx) \# \mathcal R_\si + O(d^{n-1} (\log d)^2).$$
The expectation of the number of real critical points satisfies 
now the bound 
$$ \int_{\rhl \setminus \Rr \widetilde \Delta_d} \# \text{Crit}(p_{| \Rr C_\si}) d\mu_\Rr (\si) \leq
(1 - \int_X \chi dx) \# \mathcal R_\si + O(d^{n-1} (\log d)^2).$$
Changing the function $\chi$ if necessary, the difference $(1-\int_X\chi dx) $
can be made as small as we want. We thus deduce from Proposition \ref{prop critical}
that $d^{-n} \int_{\rhl \setminus \Rr \widetilde \Delta_d} \# \text{Crit}(p_{| \Rr C_\si}) d\mu_\Rr (\si)$
converges to zero as $d$ grows to infinity and the result.

\subsubsection{Proof of Theorem \ref{theo principal produit}}
We will prove Theorem \ref{theo principal produit}
by induction on the dimension $n$ of $(X,c_X)$. 
When $n=2$, Theorem \ref{theo principal produit}
is a consequence of Theorem \ref{theo principal}.
Let us assume now that $n>2$ and that 
$(X,c_X)$ is the product of a real curve $(\Sigma, c_\Si)$ 
by a product of curves $(F,c_F)$ of dimension $n-1$. We denote by
$p : (X,c_X)\to (\Si, c_\Si)$ the projection onto the first factor. 
Again, for every section $\si \in \Rr H^0(X, L^ d)\setminus \Rr \widetilde \Delta_d$,
the restriction of $p $ to $\Rr C_\si$ satisfies the hypotheses of 
Lemma \ref{lemma Betti}, so that 
$$b_*(\Rr C_\si; \zsz) \leq 4b_*(\Rr F \cap \Rr C_\si ; \zsz) + \# \text{Crit}(p_{|\Rr C_\si}).$$
Let us bound from above each term of the right hand side of the latter
 by a $O(d^{\frac{n}{2}} (\log d)^n)$. 

Every connected component $R$ of  $\Rr X$ has a neighborhood in $X$ 
biholomorphic to a product of annuli in $\Cc$ and thus satisfies 
the conditions of Definition \ref{def asso}. Indeed, 
every complex annulus has a non vanishing holomorphic vector field.
By product, every component of $\Rr F$ has a neighborhood in $F$
trivialized by $n-1$ vector fields 
 $v_1, \cdots, v_{n-1}$. We deduce a trivialization of $\ker dp = TF$
in the neighborhood of $R$. This open neighborhood can 
be completed into an atlas adapted to $p$ with open sets disjoints from $\Rr X$.
For every $d>0$, set $\theta_d = 1-\chi_d$, where $(\chi_d)_{d>0}$ is the sequence
of functions introduced in \S \ref{para 3.3.2}. Corollary \ref{cor real expectation}
applies to $\theta_d$ whose $L^1$ norm is a $O((\frac{\log d}{\sqrt d})^n)$
whereas $||\ddbar \theta_d||_{L^1(X,\Rr)}$ is a $ O((\frac{\log d}{\sqrt d})^{n-2})$,
so that the conclusions of 
 Theorem \ref{theo equi real} also hold for $\theta_d$. As a consequence,
$$\int_{\rhl\setminus \Rr \widetilde \Delta_d} (\sum_{x \in \mathcal R_\si} \theta_d(x)) d\mu_\Rr (\si)
=O(d^{\frac{n}{2}} (\log d)^n).$$ Since for every $\si \in \rhl\setminus \Rr \widetilde \Delta_d$
and every $d>0$, 
$\# \text{Crit}(p_{|\Rr C_\si})\leq \sum_{x\in \mathcal R_\si} \theta_d(x),$
we deduce that $\# \text{Crit}(p_{|\Rr C_\si})$ is a $O(d^{\frac{n}{2}}(\log d)^n)$. 

It remains to prove that the same holds for the integral
$\int_{\rhl\setminus \Rr \widetilde \Delta_d} b_*(\Rr F \cap \Rr C_\si ; \zsz) d\mu_\Rr(\si)$. 
If the space of integration was the space of real sections of the restriction $L_{|F}$,  
this would follow from Theorem \ref{theo principal produit} in dimension $n-1$.
We will thus reduce the space of integration to this one. Let us write $(F,c_F) = (\Si_2, c_{\si_2})\times(Y,c_Y) $
where $(\Si_2, c_{\si_2})$ is a smooth real curve and $(Y,c_Y)$ a $(n-2)$-dimensional
product of curves. From Lemma \ref{lemma Betti}, for every $\si \in \rhld $,
$$b_*( \Rr F\cap \Rr C_\si; \zsz) \leq 4b_*(\Rr Y \cap \Rr C_\si ; \zsz) + \# \text{Crit}(p_{2|\Rr F\cap \Rr C_\si}),$$
where $p_2 : F\to \Si_2$ is the projection onto the first factor. Denote, with a slight abuse of notation, 
by $\theta_d$ the restriction of $\theta_d$ to $F$. We have
$$\int_{\rhld} \# \text{Crit}(p_{2|\Rr F\cap \Rr C_\si})d\mu_\Rr (\si) 
\leq \int_{\rhld} \left(\sum_{x\in \text{Crit}(p_{| F\cap C_\si})} \theta_d(x)\right)d\mu_\Rr (\si) .$$
After integration over $\rhld$ of the relation given by Proposition \ref{prop PMx} applied
to $F$ and $v_1, \cdots, v_{n-2}$, we deduce
\beq \int_{\rhld} \left(\sum_{x\in \text{Crit}(p_{| F\cap C_\si})} \theta_d(x)\right)d\mu_\Rr (\si) \leq
d^{n-1} \int_F \theta_d \omega^{n-1} +\\
  \sum_{k=0}^{n-2} (\frac{i}{2\pi})^{n-k-1}d^k
\int_F  \ddbar \theta_d \wedge \omega^k \wedge \int_{\rhld}\lambda_x (\ddbar \lambda_x)^{n-2-k} d\mu_\Rr (\si),
\eeq
 where $$\lambda_x =  
 \log \left(d^2 |f_{\si,x}|^2 + \sum^{n-2}_{i=1} |\partial f_{\si,x}(v_i)+f_{\si,x} \partial g_x(v_i)|^2\right),$$
compare with Corollary \ref{cor expectation}. We then proceed as in the proof 
of Theorem \ref{theo principal}, noting that the latter integral only depends  on the
two-jets  of  sections $\si \in \rhld$. Let us thus decompose  
$\Rhd $ as the product of the subspace $\Rr \widehat H_{3x}$ of sections
whose two-jet at $x\in F$ of their restriction to $F$ vanishes, with its orthogonal complement
$\Rr \widehat H_{3x}^\perp$. We get
$$\int_{\rhld}\lambda_x (\ddbar \lambda_x)^{n-2-k}d\mu_\Rr(\si)= 
\int_{\Rr \widehat H^\perp_{3x}\setminus \tdelta}\lambda_x (\ddbar \lambda_x)^{n-2-k}d\mu_\Rr (\si)$$
for every $k$. Likewise, the space of sections $\Rr H^0(F ; L^d_{|F})$
of the restriction $L_{|F}$ decomposes into the product $\Rr H_{3x} \times \Rr H_{3x}^\perp$.  
From Proposition \ref{prop real norm}, 
the map 
$$\si \in \Rr \widehat H_{3x}^\perp \mapsto \si_F = \frac{1}{\sqrt d} pr^\perp (\si_{|F})\in \Rr H^\perp_{3x},$$
composition of the restriction map to $F$,  a contraction by $\sqrt d$ and projection $pr^\perp $ 
onto $\Rr H^\perp_{3x}$,
 as well as its inverse are of bounded norms and determinants. 
This map in fact asymptotically coincides with the isomorphism 
$ \Rr \widehat  H^\perp_{3x} \to L^d_x\oplus T_x^* F \otimes L_x^d \oplus Sym^2(T_x^*F ) \otimes L^d_x$
given by the evaluation maps, composed with the inverse of the isomorphism 
$ \Rr H^\perp_{3x} \to   L^d_x\oplus T_x^* F \otimes L_x^d \oplus Sym^2(T_x^*F ) \otimes L^d_x$
given by the evaluation maps, see Proposition \ref{prop bounded}. From Lemma 3.1 of \cite{Tian}, 
the restricted section $\si_{|F}$ is indeed asymptotically orthogonal to $\Rr H_{3x}$.
 Hence, there is a constant $C>0$
such that 
$$ ||  \int_{\Rr \widehat H^\perp_{3x}\setminus \tdelta}\lambda_x (\ddbar \lambda_x)^{n-2-k}d\mu_\Rr (\si)      ||
\leq C  ||\int_{\Rr H^\perp_{3x}\setminus \tdelta}(\log d+ \lambda_x )(\ddbar \lambda_x)^{n-2-k}d\mu_\Rr (\si_F)||,$$
since $f_{\sqrt d \si_F,x}= \sqrt  d f_{\si_F,x}$, so that asymptotically, $\lambda_x(\si_F) = \log d + \lambda_x(\si)_{|F}$. 
Proceeding as in the proof of Theorem \ref{theo principal},
the latter right hand side is a $O(d^{n-2-k}(\log d )^3)$. 
We deduce that 
\beq
\int_{\rhld} \# \text{Crit}(p_{2|\Rr F\cap \Rr C_\si}) &=& O(d^{n-1}||\theta_d||_{L^1(F,\Rr)}) \\
&&+  O(d^{n-2}(\log d )^3||\ddbar \theta_d||_{L^1(F,\Rr)})\\
&=& O(d^{\frac{n-1}{2}}(\log d)^n)).
\eeq
Hence, the result follows by recurrence over the dimension $n$.
\subsection{Final remarks}
Several technical issues prevent us from improving Theorem \ref{theo principal}
with a $O(d^{\frac{n}{2}}(\log d)^n)$ bound in general or a $O(d^{\frac{n}{2}})$ bound.

1. First of all, Theorem \ref{theo equi real}, which 
is central in the proof of Theorem \ref{theo principal},
contains a $O(\frac{1}{d} ||\chi||_{L^1})$ term. It comes from the fact that the number 
$\# \mathcal R_\si$ of critical points of our Lefschetz pencil does
not coincide with the leading term $(\int_X c_1(L)^n)d^n$ 
given by Poincar\'e-Martinelli's formula, but is rather a polynomial 
of degree $n$ given by Proposition \ref{prop critical}
having $(\int_X c_1(L)^n)d^n$ as leading term. 
It would be of interest to identify  every monomial of the latter
with a term of Poincar\'e-Martinelli's formula or of any analytic formula.
Anyway, because of this $O(\frac{1}{d}||\chi||_{L^1})$ term in Theorem \ref{theo equi real},
we cannot use the function $\chi$ with support disjoint from $\Rr X$ which we used in the proof of 
Theorem \ref{theo principal} but rather have to use the function $\theta$ with support
in a  neighborhood of $\Rr X$ and which equals 1 on $\Rr X$ which 
we used in the proof of Theorem \ref{theo principal produit}.

2. The use of Poincar\'e-Martinelli's formula with local trivializations
$(v_1, \cdots, v_{n-1})$ forces us to choose an atlas $\mathcal U$ on $X$ and
an associated partition of unity $(\rho_U)_{U\in\mathcal U}$. As a consequence, even if
a function $\theta$ on $X$ equals one in a neighborhood of $\Rr X$,
so that $\ddbar \theta$ has support disjoint from $\Rr X$, this is not true for the functions
$\rho_U \theta$,  so that Theorem \ref{theo equi real} or Corollary \ref{cor real expectation}
cannot be used. Recall that on $\Rr X$, 
or near $\Rr X$, at
 a smaller scale than $\frac{\log d }{\sqrt d}$, the two peaks of a real peak section interfere
so that results of Proposition \ref{prop real norm} on evaluation maps no more hold. 
This forces us to have a neighborhood of each connected component of $\Rr X$ 
on which the vector fields $v_1, \cdots, v_n$ are globally defined, as 
it is the case for products of curves
 for instance.

3. Our Lefschetz pencils do not produce  Morse functions from $\Rr X$ to 
$\Rr$, but rather from $\Rr X$ to $\Rr P^1$. As a consequence, the total Betti number of 
$\Rr X$ is not bounded from above by the number of critical points of this pencil, one has to take into account also
the total Betti number of a fiber, see Lemma \ref{lemma Betti},
 and thus prove the result by induction on the dimension. Such a fiber of the pencil 
becomes then submitted to the same constraints as $X$. 

4. The weak convergence of the measure given by Theorems  \ref{theo equi}
and \ref{theo equi real} was only proved at bounded distance of the critical set of the 
Lefschetz pencil in dimensions $n>2$.  It is actually not hard
to prove it on the critical set for $n=3$ but seems less clear to us for $n>3$. This is 
another obstacle since the pencils have real critical points in general, which
have to be approached by the supports of our test functions $\theta$. 

5. The $\frac{\log d }{\sqrt d}$-scale which we use throughout the paper 
comes from Lemma \ref{lemma peak} taken out from \cite{Tian}. It ensures
that outside of the ball of radius $\frac{\log d }{\sqrt d}$, the $L^2$-norm 
of peak sections is a $O(\frac{1}{d^{2p'}})$. 
This $\frac{\log d }{\sqrt d}$ might be improved by a $\frac{1 }{\sqrt d}$ instead if a weaker
upper bound for this $L^2$-norm, such as a $O(1)$, suffices. However,
even with such a $\frac{1 }{\sqrt d}$-scale, we would still have  some $\log d$ 
term in our Theorem \ref{theo principal}, because some $\log d$
term shows up in our estimates of integrals arising from Poincar\'e-Martinelli's formula,
at the end of \S \ref{para 4.1}.\\

Note finally that  what is the exact value of the expectation in 
Theorem \ref{theo principal} remains a mystery, as well as what happens
below this expectation. Is there any exponential rarefaction 
below this expectation similar to the one  observed in \cite{GW}?
Is indeed the expectation a constant times $d^{\frac{n}{2}}$ as soon as $\Rr X$ 
is non empty?\\

\bibliographystyle{amsplain}
\bibliography{bio.bib}

\noindent
\textsc{Universit\'e de Lyon ; CNRS \\
Universit\'e Lyon 1 ; Institut Camille Jordan}\\
43 boulevard du 11 novembre 1918\\
69622 Villeurbanne cedex\\
France \\

\noindent
gayet@math.univ-lyon1.fr\\
welschinger@math.univ-lyon1.fr

\end{document}